\documentclass[10pt,a4paper]{amsart}
\usepackage{defs-eng}
\usepackage[active]{srcltx}
\usepackage{latexsym}
\usepackage{amssymb}
\usepackage{enumerate}
\usepackage[all]{xy}
\usepackage{ifpdf}

\ifpdf
        \usepackage[pdftex]{graphicx}
        \DeclareGraphicsExtensions{.pdf}
        \usepackage{hyperref}
\else
        \usepackage{graphicx}
        \DeclareGraphicsExtensions{.eps}
        \usepackage{hyperref}
\fi

\newcommand{\bA}{\mathbb A}
\newcommand{\D}{\mathcal D}
\newcommand{\C}{\mathcal C}
\newcommand{\s}{\sigma}
\newcommand{\e}{\varepsilon}
\newcommand{\Res}{\text{\rm Res}}
\newcommand{\Jac}{\text{\rm Jac}}
\newcommand{\red}{\text{\rm red}}
\newcommand{\SP}{\Sigma}

\newcommand{\nul}{{\mathfrak n}}

\title{Cohomology algebra of plane curves, weak combinatorial type, and formality}

\author[J.I.Cogolludo-Agust{\'\i}n]{Jos\'e I. Cogolludo-Agust{\'\i}n}
\address{Departamento de Matem\'aticas, IUMA\\
Universidad de Zaragoza}
\email{jicogo@unizar.es}

\author[D.Matei]{Daniel Matei}
\address{Departamento de Matem\'aticas\\
Universidad de Zaragoza \\
Institute of Mathematics "Simion Stoilow" of the Romanian Academy\\ 
P.O.Box 1-764, RO-014700 Bucharest, Romania.}
\email{daniel.matei@imar.ro}

\thanks{The first author was partially supported by the Spanish Ministry 
of Education MTM2010-21740-C02-02. The second author has been partially supported
by SB2004-0181 and grant CNCSIS PNII-IDEI 1189/2008.}

\thanks{}

\begin{document}

\begin{abstract}
We determine an explicit presentation by generators and relations of the
cohomology algebra $H^*(\PP^2\setminus \C,\CC)$ of the complement
to an algebraic curve $\C$ in the complex projective plane $\PP^2$,
via the study of log-resolution logarithmic forms on $\PP^2$.
As a first consequence, we derive that $H^*(\PP^2\setminus \C,\CC)$
depends only on the following finite pieces of data: the number of irreducible 
components of $\C$ together with their degrees and genera, the number of
local branches of each component at each singular point, and the intersection
numbers of every two distinct local branches at each singular point of $\C$.
This finite set of data is referred to as the weak combinatorial type of $\cC$.
A further corollary is that the twisted cohomology jumping loci of
$H^*(\PP^2\setminus \C,\CC)$ containing the trivial character also depend on the 
weak combinatorial type of $\cC$. Finally, the explicit construction of the generators
and relations allows us to prove that complements of plane projective curves are 
formal spaces in the sense of Sullivan.
\end{abstract}

\maketitle

\section{Introduction}
The combinatorial type $K_\C$ of a complex projective curve $\C\subset \PP^2$
consists of the following list of data: the set of irreducible components
$\C_1,\dots,\C_r$ of $\C$ together with their degrees $\bar d:=(d_1,\dots,d_r)$, 
the set of singular points $\Sing (\C)$ of $\C$ together with their topological 
types $\Sigma(\C)$, and, for every $P\in\Sing (\C)$, the correspondence $\phi_P$ 
that associates to each local branch at $P$ the global irreducible component it belongs to.

The combinatorial type of $\cC$ determines the abstract topology of $\cC$ itself. 
That is not the case for the topology of the embedding $\C\subset \PP^2$, as shown 
by Zariski's classical work where he established that the fundamental group 
$\pi_1(\PP^2\setminus \C)$ is not determined by~$K_\C$.

In this paper, we consider a topological invariant of $\cC$, called
the weak combinatorial type $W_\C$ of $\cC$, which is coarser than $K_\C$ and yet, it contains 
enough information to determine the cohomology algebra $H^*(\PP^2\setminus \C,\CC)$.
Roughly speaking, $W_\C$ consists of the following pieces of data: the set of irreducible 
components $\cC_1,\dots,\cC_r$ of $\cC$ together with their degrees $\bar d:=(d_1,\dots,d_r)$
and genera, the set of singular points $\Sing (\cC)$ of $\cC$, the correspondence $\phi_P$
as above, and the intersection numbers of every two distinct local branches at each singular 
point of $\C$. Note that the Betti numbers of $S_\cC$ depend only 
on the number, degrees, and genera of the irreducible components of~$\cC$.

In~\cite{ji-tesis}, the first author proves that $H^*(S_\cC,\CC)$ 
depends only on $W_\C$, in the case where $\cC$ is an arrangement of rational curves.
Here we extend that result to arbitrary curves.
The result follows from an explicit presentation of the cohomology algebra
$H^*(S_\C,\CC)$ which is obtained by means of the Poincar\'e residue operators
of~\cite{Griffiths-periods,Griffiths-Schmid}. For other interesting attempts to 
describe a presentation of $H^*(S_\C,\CC)$ via differential forms, see~\cite{Lubicz}.

An outline of the construction goes as follows. Fix a resolution 
$\pi: \bar{S}\to \PP^2$ of the singular locus of $\C$ in $\PP^2$ such
that $\bar{\C}$, the reduced divisor associated with $\pi^*\C$, is a simple normal crossing divisor 
in $\bar{S}$. Denote by $\bar{\C}^{[k]}$, $k=0,1,2$ the disjoint union of the codimension $k$ 
intersections of the components of $\bar{\C}$. Let $\cW_\ell^{[m]}$ be the weight filtration on the 
sheaf of logarithmic forms on $\bar{\C}^{[m]}$ with respect to the divisor $\bar{\C}^{[m+1]}$. 
Consider ${}^\ell R_m^{[k]}$ the residue operator on $\cW_\ell^{[m]}$. Note that these filtrations are
compatible with the exterior differential $d$ and that the residues ${}^\ell\tilde{R}_m^{[k]}$ defined 
on $\cW_\ell^{[m]}/\cW_{\ell-1}^{[m]}$ induce isomorphisms in $d$-cohomology. In particular, for $k=1,2$ 
one has the residue operators $R^{[k]}:={}^k R_0^{[k]}$ mapping the sheaf $\cW_k:=\cW_k^{[0]}$ of 
logarithmic forms of weight $k$ on $\bar{S}$ with respect to the divisor $\bar{\C}$ to the sheaf of 
differential forms $\cW_0^{[k]}$. Considering the complexes $(\cW_k,d)$, the exact sequence 
$0\to \cW_{i-1}\to \cW_i\to \cW_i/\cW_{i-1}\to 0$, and using the resolution $\pi$ and deRham isomorphisms, 
we can construct from the coboundary maps the following residue maps:
\[
\Res^{[i]}: H^i(S_\C,\CC)\to H^0(\bar{\C}^{[i]}), \ \ i=1,2.
\]
They are key to our approach to understand the cohomology groups $H^i(S_\C,\CC), i=1,2$.

First of all, $\Res^{[1]}: H^1(S_\C,\CC)\to H^0(\bar{\C}^{[1]})$
turns out to be an injection. Then a basis for $H^1(S_\C,\CC)$
can be chosen to be the cohomology classes of the logarithmic
$1$-forms $\sigma_i=d(\log\frac{C_i}{C_0^{d_i}}), 1\le i\le r$, where 
$\C_i$ are the irreducible components of $\C$ and $\C_0$ is a transversal
line at infinity. This condition is not strictly necessary, but we use 
it for technical reasons. For a general description see Remark~\ref{rem-affine-proj}.

The map $\Res^{[2]}: H^2(S_\C,\CC)\to H^0(\bar{\C}^{[2]})$ will not be
an injection in general, unless all the components of $\C$ are rational.
Nevertheless, we can find a decomposition of $H^2(S_\C,\CC)$
of the form $\cV^2_\cC \oplus {\mathcal K}_\cC \oplus \overline {\mathcal K}_\cC$, where
$\ker \Res^{[2]}={\mathcal K}_\C \oplus \overline {\mathcal K}_\cC$, with ${\mathcal K}_\cC$ a $g$-dimensional
vector space of classes of holomorphic $2$-forms of weight $1$ such that
$^1\tilde{R}_0^{[1]}\mathcal K_\C$ exhausts the holomorphic $1$-forms
on $\bar{\C}^{[1]}$, and $\overline {\mathcal K}_\cC$ is the conjugate of ${\mathcal K}_\cC$. 
Note that $\overline {\mathcal K}_\cC$ will necessarily consist of
classes having non-holomorphic representatives.
The vector space $\cV^2_\C$ will be generated by the classes of certain
log-resolution logarithmic $2$-forms which are constructed by the same
method employed in~\cite{ji-tesis} for the rational arrangements case.
The basic ingredients are logarithmic ideals associated with
the resolution trees appearing in the construction of $\pi$,
and ideal sheaves associated to pairs of branches at the
singular points of $\C$.
The choice of the log-resolution logarithmic $2$-forms is done by
imposing appropriate normalizing conditions.

An important feature of the decomposition
$H^2(S_\cC,\CC)=\cV^2_\cC\oplus {\mathcal K}_\cC\oplus \overline {\mathcal K}_\cC$
is that $\cV^2_\cC\supset H^1(S_\cC,\CC)\cup H^1(S_\cC,\CC)$, the cup product of
$1$-forms. Moreover, by a residue computation we determine the map
$H^1(S_\C,\CC)\times H^1(S_\C,\CC)\ \rightmap{\cup}\ \cV^2_\C$ in terms of the above 
constructed generators and see that it depends only on the degrees of the components
of $\cC$ and the intersection numbers of the local branches at the singular points 
of $\cC$. By another residue computation, we determine the
relations among the generators of $\cV^2_\cC$.
Finally, adding the trivial relations $H^1(S_\cC,\CC)\cup H^2(S_\cC,\CC)=0$,
we obtain a presentation for the cohomology algebra $H^*(S_\cC,\CC)$.

After normalizing the generators in $\cV^2_\C$ appropriately, the relations imposed by 
$H^1(S_\C,\CC)\cup H^1(S_\C,\CC)\subset \cV^2_\cC$ will be satisfied as differential forms. 
We thus derive an analogue of the Brieskorn lemma in the theory of hyperplane arrangements,
thereby obtaining an embedding of the algebra of differential forms on $S_\cC$ into the 
cohomology algebra $H^*(S_\C,\CC)$. This immediately implies the formality of $S_\C$.

The implications of this description of the cohomology algebra $H^*(S_\C,\CC)$
and the cohomology jumping loci of the space of local systems on $\CC$ 
(along the lines of~\cite{Dimca-Papadima-Suciu-formality}) will
be addressed in an upcoming paper.

\section{Settings}
\vspace*{14pt}
\subsection{$C^{\infty}$ log complex of quasi-projective algebraic varieties}
\label{sectionlogcomplex}
For the sake of completeness, we will describe in this section an appropriate 
setting for the study of the cohomology ring of the complement to plane algebraic curves. 
This includes definitions and basic properties of logarithmic sheaves and the 
definition of a very useful operator on these sheaves: the Poincar\'e residue
operator. The basic ideas in this section follow from 
\cite[Chapter 5]{Griffiths-Schmid}, but we present a slight generalization of their
results which will be necessary for the rest of the paper.

Let $X$ be a smooth, quasi-projective algebraic variety of dimension $n$ 
over $\CC$ and $\overline{X}$ be a smooth compactification of $X$. We will 
assume $\overline{X}$ to be a smooth projective variety such that 
$X=\overline{X} \setminus \overline{\D}$, where $\overline{\D}$ is a simple
normal crossing divisor, that is, a union 
$\overline{\D}_1 \cup \dots \cup \overline{\D}_N$ of smooth divisors on 
$\overline{X}$ with normal crossings. The condition of normal crossing on 
$\overline{\D}$ means that locally at $P \in \overline{X}$, the divisor 
$\overline{\D}$ is given by
$$\{(z_1,\dots,z_n) \mid z_{i_1} \cdot \ldots \cdot z_{i_m}= 0  \}= 
\{(z_1,\dots,z_n) \mid z_{I_P}=0 \},$$
where $I_P = \{i_1,\dots,i _m\} \subset \{1,\dots,n\}$. Each coordinate of $I_P$ 
must correspond locally to a different global component of $\overline{\D}$ 
(since each component $\overline{\D}_i$ is smooth). We will use a tilde, 
as in $\tilde I_P$, to indicate the ordered set of such subindices, that is, 
$\tilde I_P\subset \{1,\dots,N\}$.

\begin{definition}
\label{deflogforms}
Let ${\mathcal A}_{\overline{X}}$ be the sheaf of $C^{\infty}$ forms on 
$\overline{X}$. Denote by ${\mathcal A}^0_{\overline{X}}$ the sheaf of 
$C^{\infty}$ functions on $\overline{X}$. Note that 
${\mathcal A}_{\overline{X}}$ is a sheaf of
graded algebras over ${\mathcal A}^0_{\overline{X}}$. The 
\emph{sheaf of $C^{\infty}$ log forms}
${\mathcal A}_{\overline{X}}(\log \langle \overline{\D} \rangle)$ 
can be defined locally at a point 
$P \in \overline{X}$ as the graded algebra over
$({\mathcal A}^0_{\overline{X}})_P$ of $C^{\infty}$ forms
$\varphi \in ({\mathcal A}_X)_P$ such that
$$z_{I_P} \varphi \ \ \ {\rm and \ \ \ } z_{I_P} d\varphi$$
are $C^{\infty}$ forms in $({\mathcal A}_{\overline{X}})_P$.

A form $\varphi$ on $U\subset \overline{X}$ shall be called 
\emph{logarithmic on $U$ (with respect to $\overline{\D}$)} if 
$\varphi \in {\mathcal A}_{\overline{X}}(\log \langle \overline{\D} \rangle)(U)$.
\end{definition}

The sheaf ${\mathcal A}_{\overline{X}}(\log \langle \overline{\D} \rangle)$ 
is a locally free and finitely generated 
${\mathcal A}^0_{\overline{X}}$-algebra, as follows from

\begin{lemma}[{\rm \cite[Lemma 5.7]{Griffiths-Schmid}}]
\label{lemlogforms}
${\mathcal A}_{\overline{X}}(\log \langle \overline{\D} \rangle)(U_P) \cong 
{\mathcal A}_{\overline{X}}(U_P)\{\frac{dz_i}{z_i}\}_{i \in I_P}$
\end{lemma}
By definition,
${\mathcal A}_{\overline{X}}(\log \langle \overline{\D} \rangle)$ is 
closed under the exterior derivation $d$. This lemma shows that 
it is in fact closed under the exterior product and generated by 
${\mathcal A}^1_{\overline{X}}(\log \langle \overline{\D} \rangle)$.

In what follows, a weight filtration is defined in this sheaf 
of graded algebras that is compatible with the differential $d$.

\vspace*{10pt}
\begin{definition}
\label{defweight}
If $\ell \geq 0$ we shall define
the \emph{sheaf of $C^{\infty}$ log forms of weight $\ell$} as
the subsheaf of ${\mathcal A}_{\overline{X}}
(\log \langle \overline{\D} \rangle)$ given locally
as the $\left( {\mathcal A}_{\overline{X}} \right)^0_P$-submodule of
${\mathcal A}_{\overline{X}}(\log \langle \overline{\D} \rangle)_P$
of those forms $\varphi$ such that
$$\varphi \in 
\sum_{\smash{\mathop{I \subset I_P,}
\limits_{|I| \leq \ell}}}
{\mathcal A}_{\overline{X}} 
\left\{\frac{dz_i}{z_i} \right\}_{i \in I}.$$

\vspace*{14pt}
Such a sheaf will be denoted by 
${\mathcal W}_{\ell} :=
{\mathcal W}_{\ell} \left( {\mathcal A}_{\overline{X}}
(\log \langle \overline{\D} \rangle ) \right)$. 
If $\ell < 0$, we will assume ${\mathcal W}_{\ell} := \{0 \}$.
\end{definition}

\vspace*{10pt}
\begin{remark}
\label{remweight}
Note that 
${\mathcal W}_{\ell} \subset {\mathcal W}_{\ell+1},$
$d {\mathcal W}_{\ell} \subset {\mathcal W}_{\ell}$, and
${\mathcal W}_{\ell} \wedge {\mathcal W}_{\ell'} \subset {\mathcal W}_{\ell+\ell'}$
are obvious consequences of {\rm Definition~\ref{defweight}}.
\end{remark}

\begin{notations}
\mbox{}
\label{notfiltration}
\begin{enumerate}
\item
Let us denote by $\overline{\D}^{[k]}$ the disjoint union of the
codimension $k$ intersections of components of $\overline{\D}$, 
that is,
$$\overline{\D}^{[k]} := \bigsqcup_{|I| = k} \overline{\D}_I,$$
where $\overline{\D}_I=\cap_{i \in I} \overline{\D}_i$.
\item
There is a natural inclusion $\overline{D}_I \ \injmap{i_I} \ \overline{X}$ 
for each $\overline{\D}_I \in \overline{\D}^{[k]}$. Denoting by 
$i_k$ the corresponding map on $\overline{\D}^{[k]}$,
one has the following sheaf on $\overline{X}$
$${\mathcal A}^{*}_{\overline{\D}^{[k]}} = (i_k)_* \bigoplus_{|I|=k}
{\mathcal A}^{*}_{\overline{\D}_I}.$$
\end{enumerate}
\end{notations}

\begin{definition}
\label{defpoincresop}
Under the notations above, the 
\emph{Poincar\'e residue operator}
$$R^{[k]} : {\mathcal W}_k \left( {\mathcal A}^*_{\overline{X}}
(\log \langle \overline{\D} \rangle) \right)
\longrightarrow {\mathcal A}^{*-k}_{\overline{\D}^{[k]}}$$
can be defined locally by 
$$R^{[k]} \left( \alpha_P \wedge \frac{dz_I}{z_I} \right) =
(-1)^{\s (\tilde{I})} \alpha_P \mid_{\overline{\D}_I},$$
where:
\begin{enumerate}
\item[$i)$]
$\frac{dz_I}{z_I}$ denotes 
$\frac{dz_{i_1}}{z_{i_1}} \wedge \dots\wedge \frac{dz_{i_k}}{z_{i_k}}$, and 
\item[$ii)$]
If $\tilde I=\{\tilde i_1,\dots,\tilde i_k\}\subset \{1,\dots,N\}$, then 
$\s (\tilde{I}):=\text{sign} (\tilde i_{k+1},\dots,\tilde i_{N},\tilde i_1,\dots,\tilde i_k)$, 
where $\tilde i_{k+1} <\dots<\tilde i_{N}$ are the ordered elements of 
$\{1,\dots,N\} \setminus \tilde I$.
\end{enumerate}
\end{definition}

\begin{remark}
\label{rempoincresop}
Note that for any $\overline{\D}_{I'}$ and 
$\overline{\D}_{I}$ with $|I'| = k+1$ and $|I| = k$ one can define 
a smooth divisor on $\overline{\D}_{I}$ as follows:
$$\overline{\D}_{I} \mid_{\overline{\D}_{I'}}:=\left\{
\array{cc}
\overline{\D}_{I'} & {\rm if \ } I \subset I' \\
\emptyset & {\rm otherwise.}
\endarray \right.$$
Moreover, the union 
$$\overline{\D}_I \mid_{\overline{\D}^{[k+1]}}:=
\sum_{|I'|= k+1}
\overline{\D}_{I} \mid_{\overline{\D}_{I'}}$$
provides a simple normal crossing divisor in 
$\overline{\D}_{I} \subset \overline{\D}^{[k]}$.
Hence, $\overline{\D}^{[k]}$ can be regarded as a disjoint
union of smooth compact algebraic varieties, each component 
containing a divisor with normal crossings. Therefore, 
one can consider the sheaf of $C^{\infty}$ log forms on each 
smooth algebraic variety $\overline{\D}_{I}$ with 
respect to $\overline{\D}_I \mid_{\overline{\D}^{[k+1]}}$, denoted by 
${\mathcal A}_{\overline{\D}_I} (\log \langle \overline{\D}^{[k+1]} \rangle )$.
\end{remark}

\begin{definition}
\label{defpres2}
By means of the inclusions
$\array{ccc}
\overline{\D}_{I} & \injmap{i_k} & \overline{X}
\endarray$,
one can also define \emph{log sheaves on $\overline{\D}^{[k]}$ 
relative to $\overline{\D}^{[k+1]}$} as subsheaves of the direct 
sum of log sheaves, for each component satisfying certain 
compatibility relations. That is,
$${\mathcal A}_{\overline{\D}^{[k]}}
(\log \langle \overline{\D}^{[k+1]} \rangle) \subset
\bigoplus_{|I|= k} {\mathcal A}_{\overline{\D}_I}
(\log \langle \overline{\D}^{[k+1]} \rangle ),$$
defined by the following natural local condition:
for any strings $I_1,I_2,I'_1,I'_2$ such that
$|I_i|=|I'_i|=k_i$, $i=1,2$ and
$\{I_1 + I_2\}=\{I'_1 + I'_2\},$
and for any pair of forms
$$\alpha_P \frac{dz_{I_2}}{z_{I_2}} \in 
\left( {\mathcal A}^*_{\overline{\D}_{I_1}}
(\log \langle \overline{\D}^{[k_1+1]} \rangle) \right)_P , \quad 
\text{and} \quad 
\beta_P \frac{dz_{I'_2}}{z_{I'_2}} \in 
\left( {\mathcal A}^*_{\overline{\D}_{I'_1}}
(\log \langle \overline{\D}^{[k_1+1]} \rangle) \right)_P ,$$
one has
\begin{equation}
\label{comp-rel}
(-1)^{\sigma(\tilde{I}_1)} (-1)^{\sigma(\widetilde{I_2+I_1})}
\alpha_P \mid_{\overline{\D}_{I_1}} =
(-1)^{\sigma(\tilde I'_1)} (-1)^{\sigma(\widetilde{I'_2+I'_1})}
\beta_P \mid_{\overline{\D}_{I'_1}},
\end{equation}
where $\tilde{I}_i$ and $\tilde I'_i$ are as in Definition~\ref{defpoincresop} 
and $I+I'$ denotes juxtaposition of strings. For simplicity, this sheaf will 
be denoted by ${\mathcal A}^*_{k} (\log \langle \overline{\D} \rangle)$
and its restriction to $\overline{\D}_I$ (for $|I|=k$) by
${\mathcal A}^*_{k,I} (\log \langle \overline{\D} 
\rangle).$
\end{definition}

There also exists an obvious weight filtration 
on ${\mathcal A}_{k} (\log \langle \overline{\D} \rangle)$,
denoted by ${\mathcal W}^{[k]}_{\ell}$. Note that
${\mathcal W}^{[0]}_{\ell}={\mathcal W}_{\ell}
({\mathcal A}_{\overline{X}}(\log \langle \overline{\D} \rangle) )$
and
${\mathcal W}^{[k]}_0={\mathcal A}_{\overline{\D}^{[k]}}.$
The compatibility relations described in~{\rm (\ref{comp-rel})} 
allow for a generalization of the Poincar\'e residue operator 
to all the log sheaves relative to $\overline{\D}$, namely
\begin{equation}
\label{eq-res-gral}
{}^{\ell}R_{m}^{[k]} : {\mathcal W}^{[m]}_{\ell}
\longrightarrow {\mathcal W}^{[m+k]}_{\ell-k}.
\end{equation}
In order to give a local description of ${}^{\ell}R_{m}^{[k]}$
let us consider a point $P \in \overline{X}$ and a form 
$\varphi \in \left( 
{\mathcal A}^*_k(\log \langle \overline{\D} \rangle) \right)_P .$
Let us denote by 
$\left( {}^{\ell}R_{m}^{[k]} \varphi \right)_{I}$
the coordinate of
${}^{\ell}R_{m}^{[k]} \varphi$ on
${\mathcal A}^{*-k}_{k+m,I}(\log \langle \overline{\D} \rangle)_P$,
where $|I|=m+k$. In order to calculate this coordinate take 
two disjoint strings $I_1$ and $I_2$ such that $|I_1|=m$, and
$\overline{\D}_I=\{z_{I_1}z_{I_2}=0\}$ --\,and
hence $|I_2|=k$. The form $\varphi$ can be written as
$$\alpha \frac{dz_{I_2}}{z_{I_2}} \in \left(
{\mathcal A}^*_{k,I_1}(\log \langle \overline{\D} \rangle) \right)_P.$$
Thus one can define
$$\left( {}^{\ell}R_{m}^{[k]} \varphi \right)_{I}:=
(-1)^{\s (\tilde{I}_1)} (-1)^{\s (\widetilde{I_2+I_1})} 
\alpha \mid_{\overline{\D}_{I}}.$$
Again by~{\rm (\ref{comp-rel})} the definition of 
$\left( {}^{\ell}R_{m}^{[k]} \varphi \right)_{I}$ does not depend on
the choice of $I_1$ and $I_2$.

The main result about these generalized residue maps, which will be intensively 
used throughout the paper, is the following:

\begin{theorem}[{\rm \cite[Theorem 5.15]{Griffiths-Schmid} and 
\cite[Theorem 1.28]{ji-tesis}}]
\label{thpoincisocohomext}
Any generalized residue mapping 
$${}^{\ell}\tilde{R}^{[k]}_m : (W^{[m],*}_{\ell} / W^{[m],*}_{\ell-1}) 
\longrightarrow (W^{[m+k],*-k}_{\ell-k} / W^{[m+k],*-k}_{\ell-k-1})$$
on the complex of global sections induces an isomorphism on $d$-cohomology.
Moreover, 
$${}^{\ell-k_1}\tilde{R}^{[k_2]}_{m+k_1}\circ {}^{\ell}\tilde{R}^{[k_1]}_{m}
={}^{\ell}\tilde{R}^{[k_1+k_2]}_{m}$$
\end{theorem}

\vspace*{14pt}
\subsection
{The spaces $H^k(\PP^2 \setminus \C;\CC)$ and the residue maps}
\label{sectionresiduemap}
As a general setting, let $S$ be a smooth compact surface, $\cC\subset S$ a reduced 
divisor. Let us denote by $S_\cC$ the complement of $\cC$ in $S$. Consider a resolution 
$\pi:\overline S_\cC \to S$ of $\cC$ in $S$ such that $\overline S_\cC$ is a 
compactification of $S_\cC$ by a simple normal crossing divisor, and let $\overline \cC$ be the 
reduced structure of $\pi^*\cC$. Note that $S_\cC$ is isomorphic to 
$\overline S_\cC \setminus \overline \cC$ via $\pi$.

\begin{definition}
A \emph{log-resolution logarithmic form} on $\cC$ at $P\in S$ is a differential 
form $\varphi\in (\cA^*_{S_\C})_P$
such that $\pi^*(\varphi)_P\in \cA^*(\log \langle \overline \cC \rangle)_P$,
that is, $\varphi\in \pi_*\cA^*_{S_\C}$. The sheaf of log-resolution logarithmic forms on
$\cC$ will be denoted by $\cA_S^{\log}(\cC)=\pi_*\cA^*_{S_\cC}$.
\end{definition}

\begin{remark}
\label{rem-dominate}
Consider $\cC\subset S$ a simple normal crossing divisor, $P\in S$, and 
$\varphi\in \cA^*(\log \langle \cC \rangle)_P$ 
a differential logarithmic form. Denote by $\pi:\bar S\to S$ the blow-up of $P$ in $S$. 
Note that $\pi^*\varphi$ is also a logarithmic form on $\pi^{-1} \cC$ at any point 
$Q\in \pi^{-1}(P)$. This, together with the fact that any two sequences of blow-ups of 
$S$ are dominated by a third one, implies that the notion of log-resolution logarithmic 
form on $\cC$ is independent of the given embedded resolution of $\cC$.

Note that $\cA^{\log}_S(\cC)\subset \cA_S(\cC)$, where $\cA_S(\cC)$ is the classical 
sheaf of logarithmic differential forms on $\cC$ locally defined as 
$$(\cA_S(\cC))_P:=\{\varphi\in (\cA^*_{S_\C})_P \mid C_P \varphi\in (\cA^*_{\overline S_\cC})_P, 
\ \ C_P d\varphi\in (\cA^*_{\overline S_\cC})_P\},$$
where $C_P$ is a reduced equation of $\cC$ at $P$. 

In fact $\cA^{\log}_S(\cC)$ is the biggest subsheaf of $\cA_S(\cC)$ that is stable 
under blow-ups. Moreover, by Lemma~\ref{lemlogforms}, $\cA_S^{\log}(\cC)$ 
is locally free.
\end{remark}

\begin{construction}
\label{const-res}
Let $\C\subset \PP^2$ be an algebraic curve with irreducible components
$\C_0,\C_1,\dots,\C_r$. Fix $\pi : \overline{S}_\C \longrightarrow \PP^2$
a resolution of the singularities of $\C$ so that the reduced divisor
$\overline{\C}=(\pi^*(\C))_{\red}$ is a simple normal crossing divisor
on $\overline{S}_\C$, as described in section~{\rm \ref{sectionlogcomplex}}.

Consider the following short exact sequence of complexes
$0 \rightarrow W_{i-1} \rightarrow W_i
\rightarrow W_i/W_{i-1} \rightarrow 0$,
where $W_i$ denotes the complex
$( W^{[0],*}_i{\mathcal A}_{\overline{S}_{\cC}}
(\log \langle \overline{\cC} \rangle ),d).$ 
Let us consider its corresponding long exact sequence of $d$-cohomology
\begin{equation}
\label{grifd}
\ldots \rightarrow H^{k-1}(W_i/W_{i-1})
\rightarrow H^k(W_{i-1}) \rightarrow H^k(W_i)\
\rightmap{\delta_i^k} \ H^k(W_i/W_{i-1}) \rightarrow
H^{k+1}(W_{i-1}) \rightarrow \dots
\end{equation}
Using the de~Rham Theorem and Theorem~{\rm \ref{thpoincisocohomext}}, 
one can define the \emph{residue map} $\Res^{[i]}: H^i(S_\C;\CC)\to H^0(\overline{\C}^{[i]};\CC)$ 
as the following composition:
\begin{equation}
\label{eq-res-map}
H^i(S_\C;\CC) \congmapns{\pi^{-1}} H^i(\overline S_\C\setminus \overline \cC;\CC) \congmapns{DR}
H^i(\overline{S}_\cC,W_i) \rightmap{\delta_i^i} H^i(\overline{S}_\cC,W_i/W_{i-1})
\congmapns{R^{[i]}} H^0(\overline{\C}^{[i]};\CC)
\end{equation}
\end{construction}

\begin{proposition}[{\rm \cite[Proposition 2.2]{ji-tesis}}]
\label{prophomcohom}
If $\C$ is an algebraic plane curve with complement $S_\C$, then
$$
H^2(S_\C; \CC) \cong H_1(\C; \CC), \mbox{\ and\ } \\
H^1(S_\C; \CC) \cong H_2(\C;\CC)/\CC.
$$
\end{proposition}

\begin{notation}
Let $Y$ be a topological space. In what follows, we will
denote by $h_i(Y)$ (resp. $h^i(Y)$) the dimension of the
vector space $H_i(Y;\CC)$ (resp. $H^i(Y;\CC)$). Note that,
by the Universal Coefficient Theorem, $h_i(Y)=h^i(Y)$. 
\end{notation}

One has the following result.

\begin{proposition}
\label{propinj}
The first residue map $H^1(S_\C)\rightmap{\Res^{[1]}} H^0(\overline{\C}^{[1]})$
is injective. On the other hand
$$\ker \left(H^2(S_\C)\rightmap{\Res^{[2]}} H^0(\overline{\C}^{[2]}) \right)\subset H^2(W_1),$$
has dimension $2g$, where $g=\sum_{i=1}^r g(\C_i)$ is the sum of the genera of the 
irreducible components of $\C$.
\end{proposition}

\begin{proof}
The injectivity of $\Res^{[1]}$ follows immediately from~{\rm (\ref{grifd})} 
for the case $i=k=1$, and $H^1(W_0,d)=0$. Let us consider now~{\rm (\ref{grifd})} 
for $k=2$, $i=1$. 
\begin{equation}
\label{eq-les}
\array{ccccccc}
H^1(W_0) \to& H^1(W_1) \to & H^1(W_1/W_0) \to & H^2(W_0) \to& 
H^2(W_1) \to & H^2(W_1/W_0) \to & H^3(W_0) \\
\parallel & \parallel & \parallel & \parallel & & \parallel & \parallel\\
H^1(\overline{S}_\C) & H^1(S_\C) & H^0(\overline{\C}^{[1]}) & H^2(\overline{S}_\C) &
& H^1(\overline{\C}^{[1]}) & H^3(\overline{S}_\C) \\
\parallel & \parallel & \parallel & \parallel & & \parallel & \parallel\\
0 & \CC^r & \CC^{r+e+1} & \CC^{e+1} & & \CC^{2g} & 0 \\
\endarray
\end{equation}
where $e$ is the number of exceptional components in the resolution 
of $\C$. The equalities on the second column are a consequence of de Rham 
and Proposition~\ref{prophomcohom}. The others are a consequence of 
$H^k(W_0)=H^k(\overline{S}_\C)$, de Rham, and Theorem~\ref{thpoincisocohomext}.

Computing the Euler characteristic of this long exact sequence one obtains that
$H^2(W_1,d) \cong \CC^{2g}$ and therefore, using~{\rm (\ref{grifd})} for the case 
$i=k=2$, one obtains
$$
0 \rightarrow H^2(W_{1})=\CC^{2g} \rightarrow H^2(W_2)=H^2(S_\C)
\rightmap{\Res^{[2]}} H^0(\overline{\C}^{[2]}) \rightarrow H^{3}(W_{1}) \rightarrow \dots,
$$
which proves that 
$\ker \left( H^2(S_\C) 
\rightmap{\Res^{[2]}} H^0(\overline{\C}^{[2]}) \right)=H^2(\overline S_\cC;W_1)\cong \CC^{2g}$.
Finally, by the Leray spectral sequence, since all these sheaves are flasque, one has the projection
formula $H^2(\overline S_\cC;W_1)\cong H^2(\PP^2;\pi_* W_1)$.
\end{proof}

\subsection{Classical combinatorics}
In this paragraph we just want to give a general outline of the classical concept of \emph{combinatorial type} 
of a curve. This concept is generally accepted and used, but  is seldom explicitly defined.
In~\cite{Artal-ji-Tokunaga-survey-zariski} there is a detailed explanation of the matter. 
For the sake of completeness, we summarize the main ideas.

\begin{definition}
\label{def-comb-type}
Let $\cC\subset \PP^2$ be a plane projective curve. 
The \emph{combinatorial type} of $\mathcal C$ is given by a sextuplet
$$K_\C:=(\mathbf r, \bar d, S, \Sigma, \sigma, \Delta, \phi),$$ where:
\begin{enumerate}[$(i)$]
\item 
The elements of $\mathbf r$ are in bijection with the irreducible components of $\mathcal C$, 
\item 
$\bar d: \mathbf r \to \NN$ is the degree map that assigns to each irreducible component of 
$\mathcal C$ its degree,
\item
$S:=\Sing(\mathcal C)$, the set of singular points of $\mathcal C$, 
\item
$\Sigma$ is the set of topological types of the points in $S$,
\item
$\sigma:S\to\Sigma$ assigns to each singular point its topological type,
\item
$\Delta:=\{\Delta_P\}_{P\in S}$ where $\Delta_P$ is the set of local branches of 
$\mathcal C$ at $P\in S$, (a local branch can be seen as an arrow in the dual graph 
of the minimal resolution of $\mathcal C$ at $P$, see~\cite[Chapter~II.8]{Eisenbud-Neumann-3dimlink} 
for details) and $\phi:=\{\phi_P\}_{P\in S}$ where
$\phi_P:\Delta\to \mathbf r$ assigns to each local branch the global irreducible component 
that contains it.
\end{enumerate}
We say that two curves $\mathcal C_1$ and $\mathcal C_2$ have the 
\emph{same combinatorial type}
(or simply the \emph{same combinatorics}) if their combinatorial data $K_{\C_1})$ and $K_{\C_2}$
are equivalent, that is, if $\Sigma_{1}=\Sigma_{2},$ and there exist bijections: 
\begin{enumerate}
\item
$\varphi_{\mathbf r}:\mathbf r_1\to \mathbf r_2$,
\item 
$\varphi_{S}:S_1\to S_2$, and
\item
$\varphi_P:\Delta_{1,P}\to \Delta_{2,\varphi_{S}(P)}$ 
(the restriction of a bijection of dual graphs) for each $P\in S_1$
\end{enumerate}
such that: 
\begin{enumerate}
\item 
$\bar d_1 = \bar d_2 \circ \varphi_{\mathbf r}$, 
\item
$\sigma_{1} = \sigma_{2} \circ \varphi_{S}$, and
\item
$\varphi_{\mathbf r} \circ \phi_{1,P} = 
\phi_{2,\varphi_{S}(P)} \circ \varphi_P$.
\end{enumerate}
\end{definition}

In the irreducible case, two curves have the same combinatorial type
if they have the same degree and the same topological types for local singularities. 
On the other extreme, for line arrangements, combinatorial type is determined by the incidence graph. 
In higher dimensions, the concept of combinatorics still makes sense but it becomes much harder to 
describe, except for the case of hyperplane (or in general linear) arrangements 
where the incidence relations are enough to determined the combinatorial type.

The main interest in (and motivation for) considering combinatorial types of curves is due to the following
(see~\cite{Artal-ji-Tokunaga-survey-zariski}).

\begin{proposition}
Consider two curves $\C_1,\C_2\subset \PP^2$, and $T(\mathcal C_1)$, $T(\mathcal C_2)$ their
regular neighborhoods with boundary. Then the pairs $(T(\mathcal C_1),\C_1)$ and $(T(\mathcal C_2),\C_2)$ 
are homeomorphic if and only if $\C_1$ and $\C_2$ have the same combinatorial type.
\end{proposition}

\begin{proof}
In one direction, the self intersections of the components of $\C_i$ and the topological types of 
the singularities of $\C_i$ are well defined and preserved under homeomorphisms of pairs $(T(\mathcal C_i),\C_i)$.
This determines degrees, topological types of singularities as well as the incidence of local branches. 
Therefore their combinatorial types coincide. Conversely, the coincidence of the combinatorial
type allows one to recover the minimal resolutions of the singularities, together with a homeomorphism
between them. Since the self intersections coincide, one can extend this to a homeomorphism of the tubular 
neighborhoods of each component (including exceptional components) and glue them along the intersections 
as prescribed by the multiplicities of the components. By contracting the exceptional components one can 
define a homeomorphism of pairs between $(T(\mathcal C_1),\C_1)$ and $(T(\mathcal C_2),\C_2)$.
\end{proof}

A pair of plane curves $(\C_1,\C_2)$ such that $(T(\mathcal C_1),\C_1)$ and 
$(T(\mathcal C_2),\C_2)$ are homeomorphic, but $(\PP^2,\C_1)$ and $(\PP^2,\C_2)$ are not 
(that is, whose embeddings in $\PP^2$ are not homeomorphic) is called a \emph{Zariski pair}. 
The existence of Zariski pairs and the search for invariants of the embedding of
a curve that can tell two combinatorially-equivalent curves apart has been a very productive field of research 
started by O.Zariski in~\cite{Zariski-irregularity,Zariski-topological} (see~\cite{Artal-ji-Tokunaga-survey-zariski} 
and references therein for an extended survey on Zariski pairs).

Alternatively, one can also define a weaker concept of combinatorics as follows.

\begin{definition}
\label{def-wct}
Let $\cC\subset \PP^2$ be a plane projective curve. 
The \emph{weak combinatorial type} of $\mathcal C$ is given by a septuplet
$$W_\cC:=(\mathbf r,S,\Delta,\phi,\mu,\bar d,\bar g),$$
where
$\mathbf r$, $\bar d$, $S$, $\Delta$, and $\phi$ are defined as above,
$\bar g: \mathbf r\to \NN$ is the list of genera, and $\mu:=\{\mu_P\}_{P\in S}$, where
$\mu_{P}:\SP^2_{\phi_{P}}\Delta_{P}\to \NN$, and 
$\SP^2_{\phi_{P}}\Delta_{P}:=\frac{(\Delta_{P}\times \Delta_{P})\setminus (\Delta_{P}\times_{\phi_{P}}\Delta_{P})}{S_2}$ 
is the symmetric product of $\Delta_{P}$ outside the $\phi_{P}$-diagonal 
(that is, the fibered product $\{(\delta_1,\delta_2) \mid \phi_{P}(\delta_1)=\phi_{P}(\delta_2)\}$), 
and $\mu_{P}(\delta_1,\delta_2)$ represents the intersection multiplicity of 
$\delta_1$ and $\delta_2$ at~$P$.

We say that two curves $\mathcal C_1$ and $\mathcal C_2$ have the 
\emph{same weak combinatorial data}
(or simply the \emph{same combinatorics}) if their weak combinatorial types 
$W_{\C_1}$ and $W_{\C_2}$ are equivalent, that is, if there exist bijections: 
\begin{enumerate}
\item
$\varphi_{\mathbf r}:\mathbf r_{1}\to \mathbf r_{2}$,
\item 
$\varphi_S :S_1\to S_2$, and
\item
$\varphi_P:\Delta_{1,P}\to \Delta_{2,\varphi_{S}(P)}$ 
(restriction of a bijection of dual graphs) for each $P\in S_1$
\end{enumerate}
such that: 
\begin{enumerate}
\item 
$\bar d_{1} = \bar d_{2} \circ \varphi_{\mathbf r}$, 
\item
$\varphi_{\mathbf r} \circ \phi_{1} = \phi_{2} \circ \varphi_P$, and
\item 
$\mu_{1}(\delta_1,\delta_2)=\mu_{2}(\varphi_P(\delta_1),\varphi_P(\delta_2))$.
\end{enumerate}
\end{definition}

It is obvious that $K_\C$ determines $W_\C$ using the intersection multiplicity formula. The converse 
is also true for smooth arrangements (a curve whose irreducible components are smooth), but not true in 
general as Example~\ref{ex-wc-no-kc} shows.

The question immediately arises as to what extent the combinatorial type of a curve determines well-known
invariants of its embedding in $\PP^2$. We will refer to such invariants as \emph{combinatorial}. 
Fundamental groups of complements of curves are known not to be combinatorial, 
as first shown by Zariski in~\cite{Zariski-irregularity}. 
The cohomology ring of the complement of a curve was only known to be combinatorial when the curve 
was a line arrangement~\cite{Arnold-cohomology,Brieskorn-tresses} or more generally, 
a rational arrangement~\cite{ji-tesis}.
The purpose of the upcoming section is to prove that the cohomology ring of 
the complement of a curve is a combinatorial invariant. 

\section{Cohomology ring structure}
\label{sec-cohom-ring}
In what follows we will describe generators for $H^*(S_\C)$. For simplicity, 
we will assume $\C_0$ is a transversal line. We will consider coordinates 
$[X:Y:Z]$ in $\PP^2$ such that $\C_0:=\{Z=0\}$, and define 
$\gw:=XdY\wedge dZ+YdZ\wedge dX+ZdX\wedge dY$ the contraction of the volume form in
the affine space $\bA^3$ by the Euler vector field. 

As in the classical cases, the subspace $H^1(S_\C)$ is generated by the 
log-resolution logarithmic 1-forms $\gs_i:=d\log\frac{C_i}{C_0^{d_i}}$, 
$i=1,\dots,r$, where $C_i$ is an equation for the component $\C_i$. 

\begin{theorem}[{\cite[Theorem 2.10 and 2.11]{ji-tesis}}]
The 1-forms $\gs_i$, $i=1,\dots,r$ defined above verify the following properties:
\begin{enumerate}[$(i)$]
\item $\gs_i\in W^1 \cA^{\log}_{\PP^2}(\C)$,
\item $\gS:=\{\gs_1,\dots,\gs_r\}$ generate $H^1(S_\C)$ as a vector space, and
\item $\left(\Res^{[1]}\gs_i\right)_{\tilde\C_j}=
\begin{cases}
(-1)^{r-i} & \text{if } i=j\\
0 & \text{if } i\neq j\neq 0\\
(-1)^{r+1} d_i& \text{if } j=0.
\end{cases}$
\end{enumerate}
Moreover, $\langle \Sigma\rangle_\CC=H^1(S_\cC)=W_1^1=W^1$.
\end{theorem}

In order to obtain generators for $H^2(S_\C)$ we will define special two-forms as 
global forms of ideal sheaves. Such ideal sheaves are supported on the singular points of $\cC$.
Their definition will be given in terms of the logarithmic trees, which are isomorphic 
(as directed trees) to multiplicity trees, but whose weights are different.

\subsection{Logarithmic trees. Local setting}
Let us first recall the concept of multiplicity trees.

Let $f \in \CC\{x,y\}$ be a germ of a holomorphic function at $P$
whose set of zeroes is a curve germ $V_f \subset S_0$ with an 
isolated singularity at the point $P$. Consider the sequence of blow-ups
$$S_0 \ \longleftmap{\e_1} \ S_1 \ \longleftmap{\e_2} \ S_2 \
\longleftmap{\e_3}\ \dots\ \longleftmap{\e_m}\ S_m=\overline{S}$$
in the resolution of $S_0$ at $P$, and denote by $\pi_k$
the composition of the first $k$ blow-ups
$\pi_k=\e_k \circ \dots \circ \e_1$. The curve germ 
$\tilde{V}_{f,k}=\overline{\pi_k^{-1}(V_f \setminus \{P\})}$
shall be called the \emph{strict transform of $V_f$ in $S_k$}
and its equation denoted by $\tilde{f}_k$. 
The divisor $\pi_k^{*}(V_f)$ shall be denoted
by $\overline{V}_{f,k}$ and called the 
\emph{total transform of $V_f$ in $S_k$}. 
For simplicity, let us write $\tilde{V}_f:=\tilde{V}_{f,m}$ and 
$\overline{V}_f:=\overline{V}_{f,m}$. The exceptional divisor 
in $S_k$ resulting from the blow-up of a point in $S_{k-1}$ 
shall be denoted by $E_k$ and the points
$P_k^1,\dots,P_k^{N_k}$ in $E_k \cap \tilde{V}_{f,k}$ called
the \emph{infinitely near points to $P$ in $E_k$}. For convenience,
the point $P$ is also considered to be infinitely near to itself.
Finally, the multiplicity of $\tilde{V}_{f,k} \subset S_k$ at 
the point $P_k^\ell$ shall be denoted by $\nu_{P_k^\ell}(f)$, i.e.
$$\nu_{P_k^\ell}(f):=\mult_{P_k^\ell}(\tilde{V}_{f,k}).$$

To each resolution of singularities $\pi$ one can assign the
\emph{multiplicity tree of $\pi$ at $P$} --\,denoted by
${\mathcal T}_P(f,\pi)$, or simply by ${\mathcal T}_P(f)$ if the
resolution $\pi$ of $S_0$ is fixed. ${\mathcal T}_P(f)$ is a tree with 
weights at each vertex and is defined as follows.
\begin{enumerate}[{\rm (a)}]
\item
The vertices of ${\mathcal T}_P(f)$ are in bijection with the infinitely 
near points to $P$.
\item
Two vertices of ${\mathcal T}_P(f)$, say $Q_1$ and $Q_2$, are joined by an 
edge if and only if one of the points, say $Q_2$, belongs to $S_k$ for 
some $k$, the other point $Q_1$ belongs to $S_{k-1}$ and $Q_2 \in \e_k^{-1}(Q_1)$.
\item
For convenience, this tree is considered to simply be a vertex if $P$
is not a singular point of $f$. If $f(P) \not = 0$, then
${\mathcal T}_P(f):= \emptyset$.
\item
The weight $w({\mathcal T}_P(f),Q)$ of a vertex $Q$ is $\nu_Q(f)$. 
\end{enumerate}

\begin{example}
\label{exam-multtree}
Let $f=(x^3-y^5)(x-y^2)(y^2-x^3)y$. The tree given in Figure~{\rm \ref{figi}} 
is the multiplicity tree ${\mathcal T}_0(f)$ of the minimal resolution of $\{f=0\}$ at $0$. 

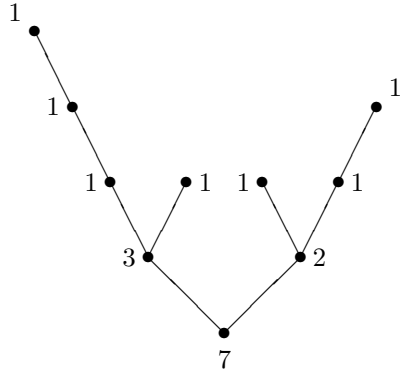
\begin{figure}[ht]
\begin{center}
\unitlength=.5mm
\begin{picture}(100,110)
\put(50,20){\circle*{3}}
\put(30,40){\circle*{3}}
\put(70,40){\circle*{3}}
\put(20,60){\circle*{3}}
\put(40,60){\circle*{3}}
\put(60,60){\circle*{3}}
\put(80,60){\circle*{3}}
\put(10,80){\circle*{3}}
\put(90,80){\circle*{3}}
\put(0,100){\circle*{3}}

\put(50,20){\line(-1,1){20}}
\put(50,20){\line(1,1){20}}
\put(30,40){\line(-1,2){10}}
\put(30,40){\line(1,2){10}}
\put(70,40){\line(-1,2){10}}
\put(70,40){\line(1,2){10}}
\put(20,60){\line(-1,2){10}}
\put(80,60){\line(1,2){10}}
\put(10,80){\line(-1,2){10}}

\put(50,13){\makebox(0,0){$7$}}
\put(25,40){\makebox(0,0){$3$}}
\put(75,40){\makebox(0,0){$2$}}
\put(15,60){\makebox(0,0){$1$}}
\put(45,60){\makebox(0,0){$1$}}
\put(55,60){\makebox(0,0){$1$}}
\put(85,60){\makebox(0,0){$1$}}
\put(5,80){\makebox(0,0){$1$}}
\put(95,85){\makebox(0,0){$1$}}
\put(-5,105){\makebox(0,0){$1$}}

\end{picture}
\end{center}
\caption{Multiplicity tree of $(f,0)$}
\label{figi}
\end{figure}
\end{example}

The set of vertices $|\mathcal T_P(f)|$ of a multiplicity 
tree $\mathcal T_P(f)$ is endowed with a partial order as follows. 
Consider $P$ as the root of the tree and direct the edges of the tree towards 
$P$. In this \emph{directed tree}, a point $Q_2$ is said to be \emph{greater} 
than $Q_1$ --\,denoted $Q_2 \geq Q_1$\,-- if there is a directed path from 
$Q_2$ to $Q_1$. In graph theory this situation is commonly described 
by calling $Q_2$ an \emph{ancestor} of $Q_1$, or $Q_1$ a \emph{descendant} 
of $Q_2$. Given a set of points $\{P_1,\dots,P_n\} \subset 
\mathcal T_P(f)$ one can define
$$\Asc(P_1,\dots,P_n)=\{Q \in {\mathcal T}_P(f) \mid 
Q \geq P_i \ \ i=1,\dots,n\},$$
and
$$\Desc(P_1,\dots,P_n)=\{Q \in {\mathcal T}_P(f) \mid 
Q \leq P_i \ \ i=1,\dots,n\}.$$
Multiplicity trees are \emph{quasi-strongly connected trees},
which means that the set of common descendants
$\Desc(P_1,\dots,P_n)$ is non empty and inherits a linear order
from ${\mathcal T}_P(f)$. The maximal element in $\Desc(P_1,\dots,P_n)$ 
is called the \emph{greatest common descendant} and is denoted by
$\gcd(P_1,\dots,P_n)$.

The \emph{degree} of a weighted tree ${\mathcal T}$ shall be defined as
\begin{equation}
\label{defdegree}
\deg({\mathcal T}):=\sum_{Q \in |{\mathcal T}| } {w({\mathcal T},Q)+1 \choose 2},
\end{equation}
where $w({\mathcal T},Q)$ denotes the weight of ${\mathcal T}$ at $Q$.

Note that if ${\mathcal T}={\mathcal T}_P(f)$, then $\deg({\mathcal T})$ is 
the $\delta$-invariant of the singularity of $f$ at $P$.

In order to simplify, we shall write ${\mathcal T} \cong {\mathcal T}'$ 
for two weighted trees that are isomorphic as trees, and 
${\mathcal T} = {\mathcal T}'$ (resp. $\geq$, $\leq$, $<$ or $>$) if 
${\mathcal T} \cong {\mathcal T}'$ and 
$\hat w({\mathcal T},Q) = \hat w({\mathcal T}',Q)$
(resp. $\geq$, $\leq$, $<$ or $>$) for any 
$Q \in |{\mathcal T}| = |{\mathcal T}'|$, where
$\hat w({\mathcal T},Q):=\sum_{Q'\in \Desc(Q)} w({\mathcal T},Q')$
(we are using the isomorphism of trees to 
identify the vertices). Note that $\hat w({\mathcal T},Q)$ is the 
multiplicity of the total transform of $f$ at $Q$. Also, ${\mathcal T} - k$ 
will denote a tree ${\mathcal T}' \cong {\mathcal T}$ so that 
$w({\mathcal T}',Q) = \max\{w({\mathcal T},Q)-k,0\}$ for 
any $Q \in |{\mathcal T}|$. 
Particularly useful will be the tree
\begin{equation}
\label{defT0}
{\mathcal T}^{\nul}_P(f) := {\mathcal T}_P(f) - 1.
\end{equation}
Sometimes it will be necessary to compare empty trees.
In this case, the conditions
$=, \leq, \geq$ are vacuous and hence always satisfied.

Let $g \in \CC\{x,y\}$ be another germ at $P$.
Then one can consider the \emph{restriction of $g$ to a weighted tree $\mathcal T$ 
(e.g. ${\mathcal T}_P(f)$}) --\,denoted by ${\mathcal T}|_g$\,-- as a
weighted graph isomorphic to ${\mathcal T}$, whose weight at
each vertex $Q$ is $\nu_Q(g)$. One can check that the set 
$I:=\{g \in \CC\{x,y\} \mid {\mathcal T}|_g \geq {\mathcal T}\}$
defines an ideal. Note that,
$\dim_{\CC} \left( \CC\{x,y\}/{\mathfrak m}^k \right)={k+1 \choose 2}$.
Hence
\begin{equation}
\label{eq-deg-ideal}
\deg({\mathcal T}) = \dim_{\CC} \frac{\CC\{x,y\}}{ I }.
\end{equation}

Let $\C$ be a plane projective curve and $P\in\Sing\C$. Note that its multiplicity tree 
does not depend on the equation of $\D$, hence it can be denoted as $\mathcal T_P(\C)$ 
or $\mathcal T_P(\C,\pi)$ in case we want to specify the underlying resolution.

In case $\C$ is irreducible and has degree $d$, 
from~{\rm (\ref{defdegree})} and~{\rm (\ref{defT0})} one 
can rewrite the Noether formula for the genus 
\cite[p.\,614]{Brieskorn-Knorrer-curves} as follows
\begin{equation}
\label{genusform}
g(\C)= \frac{(d-1)(d-2)}{2} - \sum_{P \in \Sing \C}
\deg({\mathcal T}^{\nul}_P(\C)).
\end{equation}

Now consider $\pi$ a resolution of singularities for the plane curve $\C$. 
We define the \emph{basic ideal sheaf of $\C$ with respect to $\pi$} as follows
\begin{equation}
\label{basic-ideal-sheaf}
({\mathcal I}^{\nul}_{\C,\pi})_P:=\{ h \in {\mathcal O}_P \mid 
{\mathcal T}_P(\C,\pi) |_h \ \geq \ {\mathcal T}^{\nul}_P(\C,\pi)\}.
\end{equation}
If no possible confusion results from the underlying 
resolution, the sheaf ${\mathcal I}^{\nul}_{\C,\pi}$ will be 
denoted simply by ${\mathcal I}^{\nul}_{\C}$.

\begin{remark}
Since $\pi$ also induces a log resolution of the ideal $C=I(\C)$ at any point $P$, one can also see 
${\mathcal I}^{\nul}_{\C}$ as the multiplier ideal sheaf of $C$, that is 
$\pi_* \mathcal O_{\bar S_{\C}}(K_{\bar S_\C / \PP^2}-F)$, 
where $C \cdot \mathcal O_{\bar S_{\C}}=\mathcal O_{\bar S_{\C}}(-F)$.
Analogously, ${\mathcal I}^{\nul}_{\C}$ corresponds to the special ideal of quasi-adjunction $\mathcal A_{0}(C)$
as defined in~\cite{Libgober-alexanderinvariants,Loeser-Vaquie-Alexander}.
\end{remark}

This leads naturally to the idea of logarithmic ideal of a germ.
Let $f \in \CC \{x,y\}$ be a holomorphic germ at $P$ and $\pi$ a resolution of $V_f$. 
We can define an ideal $I \subset \CC \{x,y\}$ satisfying that for any germ $h \in I$ the 2-form
\begin{equation}
\label{local2form}
h \ \frac{dx \wedge dy}{f}
\end{equation}
is log-resolution logarithmic at $P$ --\,with respect to $V_f$ and the
resolution $\pi$.
Such an ideal will be called a \emph{logarithmic ideal for $f$ at $P$}.

\begin{remark}
\label{rem-log-ideals-indep}
Using Remark~\ref{rem-dominate}, it is easy to see that the set of logarithmic ideals associated with a 
singularity is independent of the resolution. Hence, from now on, when referring to logarithmic ideals
any reference to a resolution will be omitted.
\end{remark}

A very useful way to encode the information required to construct logarithmic ideals is given by weighted trees.

\begin{definition}
\label{deflogtree}
Let $(f,0)$ be a germ and $\mathcal T(f)$ its multiplicity tree. 
A weighted tree ${\mathcal T}$ is said to be a \emph{logarithmic tree for $(f,0)$} 
if it satisfies the following properties:
{\noindent
\begin{enumerate}[$(i)$]
\item
${\mathcal T} \cong {\mathcal T}(f),$ and
\item
the ideal
$I:=\{h \in {\mathcal O}_0 \mid 
{\mathcal T}(f)|_h \ \geq \ {\mathcal T} \}$
is logarithmic.
\end{enumerate}
}
In addition, if $\delta_1$ and $\delta_2$ are local branches of $(f,0)$, we say that 
${\mathcal T}$ is a \emph{logarithmic tree $($for $(f,0))$ relative to $\delta_1$ and $\delta_2$}
if ${\mathcal T}$ satisfies $(i)$, $(ii)$ above and
{\noindent
\begin{enumerate}
\item[$(iii)$]
\label{deflogtree-MI}
if $\varphi\in M_I$, where 
$M_I:=\{h \in {\mathcal O}_0 \mid {\mathcal T}(f)|_h \ = \ {\mathcal T} \} \subset I$, 
then
$$\left( \Res^{[2]} \varphi \
\frac{dx \wedge dy}{f} \right)_{Q} \not = 0$$
if and only if $Q$ is a vertex of the unique subtree 
$\gamma(\delta_1,\delta_2)\subset \cT$ joining $\delta_1$ and $\delta_2$.
\end{enumerate}
}
\end{definition}

\begin{example}
\label{ex-inul}
Note that $\cT^{\nul}_P(f)$ is a logarithmic tree for $f$, but it is \emph{not} relative to 
any two branches $\delta_1$, $\delta_2$. One can check that properties $(i)$ and 
$(ii)$ are satisfied, but $\cT^{\nul}_P(f)$ does not satisfy property $(iii)$. Moreover, if 
$\varphi\in \CC\{x,y\}$ is a germ at $P$ such that $\psi:=\varphi\frac{dx \wedge dy}{f}$ with 
$\cT_P(f)|_\varphi \geq \cT^{\nul}_P(f)$, then one can check that $\psi\in \pi_* W^2_1$, that 
is, it has weight one and hence $\left(\Res^{[2]}\psi\right)_Q=0$ for any vertex $Q$ of $\cT_P(f)$.
\end{example}

The main result of this part is the existence of logarithmic trees relative to any two local branches 
of any reduced germ $f$.

\begin{theorem}[{\cite[Lemma 2.34]{ji-tesis}}]
\label{thm-ex-logtree}
For any given two local branches $\delta_1$ and $\delta_2$ of $f$ at $P$,
there exists a logarithmic tree for $(f,P)$ relative to $\delta_1$ and $\delta_2$.
\end{theorem}

Theorem~\ref{thm-ex-logtree} is constructive. We will denote such a tree by $\cT_P^{\delta_1,\delta_2}$ 
and it will be referred to as the \emph{basic logarithmic tree relative to $\delta_1$ and $\delta_2$}. 

\subsection{Logarithmic ideal sheaves. Global setting.}
Let us return to the situation presented at the beginning of this section, where $\C$ is a 
plane projective curve and $\pi$ is a resolution of singularities. The concept of logarithmic ideal 
translates globally as follows:

\begin{definition}
We call an ideal sheaf ${\mathcal I}$ on $\PP^2$ \emph{logarithmic ideal sheaf} for $\C$ if its 
stalks ${\mathcal I}_P$ are logarithmic ideals for the germs $C_P$ of $\C$ at any $P \in \PP^2$.
\end{definition}

\begin{remark}
By Example~\ref{ex-inul}, the basic ideal sheaf of $\C$ denoted by ${\mathcal I}^{\nul}_{\C}$ 
defined in~(\ref{basic-ideal-sheaf}) is a logarithmic ideal sheaf for $\C$.

Also note that by Remark~\ref{rem-log-ideals-indep}, such sheaves are independent of the given resolution.
\end{remark}

Let $\C_{ij}:=\cC_i\cup \cC_j$, $d_{ij}:=\deg \cC_{ij}$, and $g_{ij}:=g(\C_{ij})$. We denote by $C_{ij}$ 
an equation of $\cC_{ij}$ (which is $C_i$ if $i=j$ and $C_iC_j$ if $i\neq j$). We will first check that
the basic ideal sheaf of $\cC_{ij}$ has non-trivial sections of degree $d_{ij}-2$ except for the obvious 
case of lines.

\begin{proposition}
\label{prop-dimnij}
 $\dim H^0(\PP^2,{\mathcal I}_{\cC_{ij}}^{\nul}(d_{ij}-2)) \geq d_{ij}+g_{ij}-\#\{i,j\}$.
\end{proposition}

\begin{proof}
To ease the notation let us write $\cI$ for $\cI_{\cC_{ij}}^{\nul}$. From the exact sequence
$$
0 \to \cI(d_{ij}-2) \to \cO_{\PP^2}(d_{ij}-2) \to \cO/\cI(d_{ij}-2) \to 0
$$
and the fact that $H^\ell(\cO(k))=0$ for any $k$ and $\ell>0$, one obtains that
\begin{equation}
\label{eq-h0-1}
h^0(\PP^2,\cI(d_{ij}-2))\geq {d_{ij} \choose 2} - h^0(\PP^2,\cO/\cI).
\end{equation}
In what follows, we will assume $i\neq j$. The case $i=j$ is analogous. First, we will calculate 
$h^0(\PP^2,\cO/\cI)$. Note that, by~(\ref{eq-deg-ideal}), one has 
$$h^0(\PP^2,\cO/\cI) = \sum_{P\in \Sing \cC_{ij}} \deg \cT^{\nul}_P(\cC_{ij})=
\sum_{P\in \Sing \cC_{ij}} \sum_{Q\in |\cT_P(\cC_{ij})|} {\nu_Q(C_{ij}) \choose 2}.
$$
Since $\nu_Q(C_{ij})=\nu_Q(C_{i})+\nu_Q(C_{j})$ and ${a+b \choose 2}={a \choose 2}+{b \choose 2}+ab$
one obtains that 
$$
h^0(\PP^2,\cO/\cI) = 
\sum_{P\in \Sing \cC_{ij}} \sum_{Q\in |\cT_P(\cC_{ij})|} {\nu_Q(C_{i}) \choose 2} + {\nu_Q(C_{j}) \choose 2}
+ \nu_Q(C_{i})\nu_Q(C_{j})
$$
and finally using~(\ref{genusform}) one obtains
\begin{equation}
\label{eq-h0-2}
h^0(\PP^2,\cO/\cI) = 
{d_i-1 \choose 2}+{d_j-1 \choose 2}-g_{ij}+d_id_j.
\end{equation}
Therefore~(\ref{eq-h0-1}) becomes
$$
h^0(\PP^2,\cI(d_{ij}-2))\geq \left[ {d_{ij} \choose 2} - {d_i-1 \choose 2} - {d_j-1 \choose 2} - d_id_j \right] 
+ g_{ij}= (d_i-1)+(d_j-1)+g_{ij}.
$$
\end{proof}

In order to construct global forms we will proceed as follows. 
First, for each irreducible component $\cC_i$ of $\cC$, we will order the $d_i=\deg \cC_i$
points of $\cC_i$ at infinity $\cC_i\cap \cC_0=\{P_1^i,\dots,P_{d_i}^i\}$. 

\begin{definition}
\label{defIdelta}
Let $P \in \C_{ij}$, and let $\delta_1$ (resp. $\delta_2$) be a local
branch of the irreducible component $\C_i$ (resp. $\C_j$) at $P$.
The \emph{ideal sheaf ${\mathcal I}^{\delta_1,\delta_2}_{\cC_{ij}}$
associated with $\delta_1$ and $\delta_2$} shall be defined as
$$({\mathcal I}_{\cC_{ij}}^{\delta_1,\delta_2})_Q:=
\left\{ h \in {\mathcal O}_Q \left|
\array{ll}
{\mathcal T}_Q(\C_{ij})|_h \ \geq \ {\mathcal T}_P^{\delta_1,\delta_2}(\C_{ij}) & 
{\rm if \ } Q = P, \\
{\mathcal T}_Q(\C_{ij})|_h \ \geq \ {\mathcal T}_Q(\C_{ij})-2 & {\rm if \ } Q \in \{P_1^i,P_1^j\},\\
{\mathcal T}_Q(\C_{ij})|_h \ \geq \ {\mathcal T}_Q^{\nul}(\C_{ij}) & {\rm otherwise \ }
\endarray
\right.
\right\}.$$
A global section $s$ of ${\mathcal I}_{\C_{ij}}^{\delta_1,\delta_2}(d)$ 
shall be called \emph{essential} if $s_Q \in M_{I_Q}$ for
every $Q \in \PP^2$, where $s_Q$ is the section $s$ localized at
$Q$, $I_Q=({\mathcal I}_{\C_{ij}}^{\delta_1,\delta_2})_Q$ and
$M_{I_Q}$ is as in Definition~{\rm \ref{deflogtree}}.
\end{definition}

Analogously to~\cite[Lemma 3.35]{ji-tesis} one can prove the following.

\begin{proposition}
\label{prop-degij}
 $\deg {\mathcal I}_{\C_{ij}}^{\delta_1,\delta_2}=
\deg {\mathcal I}_{\C_{ij}}^{\nul}+d_{ij}-\#\{i,j\}-1$.
\end{proposition}

Therefore Propositions~\ref{prop-dimnij} and~\ref{prop-degij} imply the following.

\begin{proposition}
\label{prop-dimij}
 $\dim H^0(\PP^2,{\mathcal I}_{\C_{ij}}^{\delta_1,\delta_2}(d_{ij}-2))>g_{ij}$.
\end{proposition}

One can give a description of a section in such sheaf ideals. 
In order to do so, let us denote by $\gamma_P(\delta_1,\delta_2)$ the minimal 
subtree in $\cT_P(\C_{ij})$ (see Definition~\ref{deflogtree}) containing $\delta_1$ 
and $\delta_2$. We can consider $\gamma_P(\delta_1,\delta_2)$
as a subset of the total transform $\overline \C$ of $\C$ (in particular it should 
contain $\tilde \C_i$ and $\tilde \C_j$). We also denote by 
$v(\gamma_P(\delta_1,\delta_2))$ the set of vertices of $\gamma_P(\delta_1,\delta_2)$.

\begin{proposition}
\label{prop-resij}
Let $\varphi$ be a section in
$H^0(\PP^2,{\mathcal I}_{\C_{ij}}^{\delta_1,\delta_2}(d_{ij}-2))$. Consider the 
2-form $\varphi \frac{\omega}{C_0C_{ij}}$. One has the following basic properties:
\begin{enumerate}[$(1)$]
\item \label{prop-resij-res}
$$\left(\Res^{[2]}\varphi \frac{\omega}{C_0C_{ij}} \right)_Q=
\begin{cases}
\pm \lambda & \text{if } Q\in v(\gamma_P(\delta_1,\delta_2)),\\
\pm \gve_{ij}\lambda & \text{if } Q\in \{P_1^i,P_1^j\},\\
0 & \text{otherwise,}
\end{cases}
$$
where 
$
\gve_{ij}=\begin{cases}
1 & \text{if } i\neq j \\
0 & \text{otherwise.}
          \end{cases}
$

Moreover, $\lambda\neq 0$ if and only if 
$\varphi\in M_{{\mathcal I}_{C_{ij}}^{\delta_1,\delta_2}}$ is essential 
(as in~Definition~{\rm \ref{deflogtree}}). 
\item \label{prop-resij-signs}
The signs of the residues described in~$(\ref{prop-resij-res})$ are such that if 
$\D\subset \gamma_P(\delta_1,\delta_2) \cup \C_0 \subset \overline{\C}^{[1]}$ is an irreducible
component then
$$
\left({}^2\tilde{R}^{[1]}_0 \varphi \frac{\omega}{C_0C_{ij}}\right)_F=
\begin{cases}
0 & \text{ if } F \notin \D \\
\lambda \omega_1 & \text{ if } F \in \D, F\neq \tilde \C_0, \tilde \C_i, \tilde \C_j\\
\gve_{ij} \lambda \omega_1 & \text{ if } F = \tilde \C_0 \\
\eta_i^{\delta_1,\delta_2} & \text{ if } F = \tilde \C_i,\\
\eta_j^{\delta_1,\delta_2} & \text{ if } F = \tilde \C_j,\\
\end{cases}
$$
where $\omega_1:=\frac{XdY+ YdX}{XY}$ is the Euler 1-form on $\PP^1\cong F$ that
has poles on $[0:1], [1:0]$ and whose residues are $\lambda$ and $-\lambda$ respectively,
and where $\eta_i, \eta_j$ are 1-forms on $\tilde \C_i$ (resp. $\tilde \C_j$) with only two poles.
Moreover, the poles are at the point on $\tilde \cC_i$ (resp. $\tilde \cC_j$) determined by $\delta_1$ 
(resp. $\delta_2$) and at the point determined by the branch of $\tilde \cC_i$ (resp. $\tilde \cC_j$) 
at $P_1^i$ (resp. $P_1^j$).
\end{enumerate}
\end{proposition}

\begin{proof}
Since the proof of Theorem~\ref{thm-ex-logtree} is constructive, one can use such a construction
to check part~(\ref{prop-resij-res}). Part~(\ref{prop-resij-signs})
is a consequence of the commutativity of the generalized residue maps
(Theorem~\ref{thpoincisocohomext}), and the fact that the residues of a meromorphic
function on a compact Riemann surface add up to zero and the difference between the number of
zeroes and poles is the Euler characteristic $2g-2$, where $g$ is the genus.
\end{proof}

By Theorem~\ref{thpoincisocohomext} and the exact sequence~(\ref{eq-les}) note that 
\begin{equation}
\label{eq-sec}
H^2(W_1) \cong H^2(W_1/W_0) \congmap{{}^{1}\tilde{R}^{[1]}_0}
H^1(W_{0}^{[1]})=H^1(\overline{\C}^{[1]}).
\end{equation}
Using the inclusion $\Omega^*\ \injmap{i^*}\ W^*_{0}$, from the complex of global holomorphic 
forms on $\overline{\C}^{[1]}$ to the complex of global differential forms, one has a map 
$H^1({\overline{\C}^{[1]}},\Omega^1) \ \rightmap{i^1}\ H^1(\overline{\C}^{[1]})$.
Also note that $\dim H^1({\overline{\C}^{[1]}},\Omega^1)=g$. In the following, we will
describe generators for $H^1({\overline{\C}^{[1]}},\Omega^1)$.

\begin{proposition}
\label{prop-ki}
Let  
$K_i:=\{\psi=\varphi\frac{\omega}{C_i} \mid \varphi\in H^0(\PP^2,\mathcal O(d_i-3)), 
\cT_P(\C_i)|_\varphi\geq \cT^{\nul}_P(\C_i)\}$ and $K_\cC=\sum_{i=1}^r K_i$.
One has the following properties:
\begin{enumerate}[$(1)$]
 \item \label{prop-ki-log}
$K_i\subset W^2_1\left( \cA^{\log}_{\PP^2}(\C) \right)$,
 \item \label{prop-ki-ker}
$\mathcal K_\C\subset \ker \left( H^2(S_\C) 
\rightmap{\Res^{[2]}} H^0(\overline{\C}^{[2]}) \right)$, where 
$\mathcal K_\cC$ is the projection of
$K_\cC$ on $H^2(W_1)\subset H^2(S_\C)$,
 \item \label{prop-ki-oplus}
$K_\cC=\oplus_{i=1}^r K_i$
and 
$\cK_\cC=\oplus_{i=1}^r \cK_i$, 
 \item \label{prop-ki-im}
${{}^{1}\tilde{R}^{[1]}_0} \mathcal K_\C = i^1 H^1({\overline{\C}^{[1]}},\Omega^1).$
\end{enumerate}
Moreover, if $\overline{K}_\cC$ denotes the conjugate of $K_\cC$, then 
$\mathcal K_\C \oplus \overline{\mathcal K}_\C = 
\ker \left( H^2(S_\C) \rightmap{\Res^{[2]}} H^0(\overline{\C}^{[2]}) \right)$.
\end{proposition}

\begin{proof}
Let us start with parts~(\ref{prop-ki-log}) and~(\ref{prop-ki-ker}). 
The result is local, and according to 
Example~\ref{ex-inul}, it is enough to check it at the points at infinity 
$\{P^i_1,\dots,P^i_{d_i}\}=\C_0\cap \C_i$. Since such points are smooth on $\C_i$
by hypothesis, the condition $\cT_{P^i_k}(\C_i)|_\varphi\geq \cT^{\nul}_{P^i_k}(\C_i)$
is vacuous. Hence the local equation of $\psi$ at ${P^i_k}$ is 
$$
\varphi(u,v)\frac{du\wedge dv}{v},
$$
up to a unit, where $\{v=0\}$ (resp. $\{u=0\}$) is the local equation 
of $\C_i$ (resp. of $\C_0$). Hence part~(\ref{prop-ki-log}) follows, as well as 
the first statement of part~(\ref{prop-ki-ker}).
On the one hand note that 
$$\dim K_i\geq {{d_i-1}\choose{2}} - \sum_{P\in \Sing(\C_i)}\deg \cT_P^{\nul}(\C_i)=g_i.$$
In order to prove the first statement of part~(\ref{prop-ki-oplus}), note that if
$\psi_1+\dots+\psi_r=0$, $\psi_i\in K_i$, then multiplying by $C$, one has 
\begin{equation}
\label{eq-c}
\varphi_1C_2\cdot\ldots\cdot C_r+C_1\varphi_2\cdot \ldots 
\cdot C_r+\dots+C_1C_2\cdot \ldots \cdot \varphi_r=0.
\end{equation} 
Consider 
$\{Q_1,\dots,Q_{d_1-2}\}\subset \C_1\setminus (\C_2\cup \dots \cup\C_r)$ and evaluate such
points in~(\ref{eq-c}). One obtains that $\varphi_1(Q_1)=\dots=\varphi_1(Q_{d_1-2})=0$.
Since $\C_1$ is irreducible and $\deg \varphi_1=d_1-3$, one obtains that $\varphi_1=0$.
Proceeding analogously for every $\varphi_i$, one obtains that 
$K_\C=K_1\oplus K_2\oplus \dots \oplus K_r$. This same idea shows that 
${{}^{1}\tilde{R}^{[1]}_0} \psi=0$ if and only if $\psi=0$. Therefore
$\dim {{}^{1}\tilde{R}^{[1]}_0} \mathcal K_\C \geq g = h^1({\overline{\C}^{[1]}},\Omega^1)$.
The inclusion ${{}^{1}\tilde{R}^{[1]}_0}\mathcal K_\C\subset 
i^1 H^1({\overline{\C}^{[1]}},\Omega^1)$
forces $i^1$ to be an injection and thus, the second statement of part~(\ref{prop-ki-oplus})
and part~(\ref{prop-ki-im}) follow for dimension reasons.

The \emph{moreover} statement is a consequence of Proposition~\ref{propinj}.
\end{proof}

\begin{remark}
Note that Proposition~\ref{prop-ki}(\ref{prop-ki-im}) implies in particular that 
cohomology classes outside $\mathcal K_\C$ do not have holomorphic representatives.
\end{remark}

\begin{notation}
\label{not-vaff}
One can normalize any log-resolution logarithmic 2-form 
$\varphi \frac{\omega}{C_0C_{ij}}$ as in Proposition~\ref{prop-resij} 
in such a way that $\left(\Res^{[2]}\varphi \frac{\omega}{C_0C_{ij}}\right)_{P^i_1}=1$.

Note that if $\psi_{ij}\in K_{ij}=K_i \oplus K_j$, then
$\left(\Res^{[2]}\varphi \frac{\omega}{C_0C_{ij}}+\psi_{ij}\right)_{P^i_1}=1$.
The set of classes of such normalized 2-forms will be denoted by 
$v_1:=\{\psi_P^{\delta_1,\delta_2}+K_{ij}\}_{P,\delta_1,\delta_2}$.
\end{notation}

Analogously, one needs to consider certain 2-forms with residues on the line at infinity.

\begin{definition}
\label{defIinfty}
For any $P_k^i\in \C_0\cap \C_i$, $k=2,\dots,d_i$,
the \emph{ideal sheaf ${\mathcal I}^{P_k^i}_{\C_{i}}$
associated with $P_k^i$} shall be defined as
$$({\mathcal I}_{\C_{i}}^{P_k^i})_Q:=
\left\{ h \in {\mathcal O}_Q \mid 
\array{ll}
{\mathcal T}_Q(\C_{i})|_h \ \geq \ {\mathcal T}_Q(\C_{i})-2 & {\rm if \ } Q \in \{P_1^i,P_k^i\}\\
{\mathcal T}_Q(\C_{i})|_h \ \geq \ {\mathcal T}_Q^{\nul}(\C_{i}) & {\rm otherwise \ }
\endarray
\right\}.$$
As above, a global section $s$ of ${\mathcal I}_{\C_{i}}^{P_k^i}(d)$ 
shall be called \emph{essential} if 
${\mathcal T}_Q(\C_{i})|_{s_Q}={\mathcal T}_Q(\C_{i})-2$ for
every $Q \in \{P_1^i,P_k^i\}$, where $s_Q$ is the section $s$ localized at $Q$.
\end{definition}

One can also describe such sections in terms of their residues as in 
Proposition~\ref{prop-resij}. Its analogue reads as follows.

\begin{proposition}
\label{prop-resinfty}
Let $\varphi$ be a section in
$H^0(\PP^2,{\mathcal I}_{\C_{i}}^{P_k^i}(d_{i}-2))$. Consider the 
2-form $\psi_\infty^{i,k}:=\varphi \frac{\omega}{C_0C_{i}}$.
One has the following basic properties:
\begin{enumerate}[$(1)$]
 \item \label{prop-resinfty-res}
$$\left(\Res^{[2]}\psi_\infty^{i,k}\right)_Q=
\begin{cases}
\pm \lambda & \text{if } Q\in \{P_1^i,P_k^i\},\\
0 & \text{otherwise.}
\end{cases}
$$
 \item \label{prop-resinfty-l}
$\lambda\neq 0$ if and only if $\varphi\in M_{{\mathcal I}_{\C_{i}}^{P_k^i}}$ 
is essential (as defined in~{\rm \ref{defIinfty}}). 
\item \label{prop-resinfty-signs}
The signs of the residues described in~$(\ref{prop-resinfty-res})$ are such that if 
$\D\in \{\C_i,\C_0\} \subset \overline{\C}^{[1]}$, 
then ${}^2\tilde{R}^{[1]}_0 \psi_\infty^{i,k}$ has exactly two poles along 
$\D$ whose residues are $\lambda$ and $-\lambda$ so that they add up to zero.
\end{enumerate}
\end{proposition}
\begin{proof}
The proof follows immediately from~\ref{ex-inul} and the fact that the singularities
at infinity are always nodes, hence the local trees are as shown in Figure~\ref{fig-node}.
\begin{figure}[ht]
\includegraphics[scale=.5]{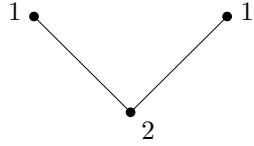}
\begin{picture}(0,0)
\put(-35,-5){\makebox(0,0){$2$}}
\put(-85,40){\makebox(0,0){$1$}}
\put(2,40){\makebox(0,0){$1$}}
\end{picture}
\caption{Multiplicity tree of a node.}
\label{fig-node}
\end{figure}
Therefore $\cT_{P^i_1}(\C_i)-2$, imposes no conditions on $\varphi$ and thus
$$
\left(\Res^{[2]}\psi_\infty^{i,k}\right)_{P^i_1}=\pm \varphi(P^i_1)=\pm \lambda.
$$
The same argument works for $P^i_k$, and 
hence~(\ref{prop-resinfty-res})-(\ref{prop-resinfty-l}) follow. 
Finally, part~(\ref{prop-resinfty-signs}) follows from the same ideas as 
in Proposition~\ref{prop-resij}(\ref{prop-resij-signs}).
\end{proof}

\begin{notation}
\label{not-vinfty}
For any $P^i_k\in \C_0\cap \C_i$, $k=2,\dots,d_i$, let 
$\psi_\infty^{i,k}$ be a log-resolution logarithmic 2-form as constructed in 
Proposition~\ref{prop-resinfty}, again, with the extra normalizing condition that 
$\left(\Res^{[2]}\psi_\infty^{i,k}\right)_{P^i_1}=1$. Again, as in Notation~\ref{not-vaff},
the set of all classes of such normalized 2-forms will be denoted by 
$v_\infty:=\{\psi_\infty^{i,k}+K_{i}\}_{i,k}$.
\end{notation}

Note that in the line arrangement case, $v_\infty$ is trivial and each class in $v_1$ contains 
exactly one representative, which is the form $\frac{d\ell_i\wedge d\ell_j}{\ell_i\ell_j}$ 
described by Brieskorn~\cite{Brieskorn-tresses} and Orlik-Solomon~\cite{Orlik-Solomon}. 
Note that their cohomology classes are not linearly independent. In particular, 
if $\ell_i,\ell_j$, and $\ell_k$ intersect at $P$, then the 2-form
\begin{equation}
\label{eq-lines}
\frac{d\ell_i\wedge d\ell_j}{\ell_i\ell_j}+\frac{d\ell_j\wedge d\ell_k}{\ell_j\ell_k}+
\frac{d\ell_k\wedge d\ell_i}{\ell_k\ell_i}
\end{equation}
is not only cohomologically trivial, but even 0 as a 2-form as one can easily see as follows:
by the concurrence condition we can assume $\ell_i=x$, $\ell_j=y$, and $\ell_k=\alpha x+\beta y$.
Hence $d\ell_i\wedge d\ell_j=dx\wedge dy$, $d\ell_j\wedge d\ell_k=-\alpha dx\wedge dy$, and
$d\ell_k\wedge d\ell_i=-\beta$. Thus, after removing denominators~(\ref{eq-lines}) becomes
$(\ell_k-\alpha x-\beta y)dx\wedge dy$ which is zero.

The following results will be used as a tool to come up with the relations amongst 
cohomology classes. Additionally, we will prove that such relations also hold for the 
2-forms in $v_1$ and not just for their cohomology classes.

\begin{lemma}
\label{lema-key}
Let $\psi:=\varphi\frac{\gw}{q C_i}\in W^2_1$. If $\left({}^1\tilde R^{[1]}_0(\psi)\right)_{\tilde \C_i}=0$, 
then $\varphi=p C_i$ for some function $p$. Moreover, if $\varphi$, $q$, and $C_i$ are homogeneous 
polynomials, then $p\in \CC[X,Y,Z]$ is a homogeneous polynomial.
\end{lemma}

\begin{proof}
Note that for any $P\in \C_i\setminus \Sing \C$, one has that 
$\left({}^1R^{[1]}_0(\psi)\right)_{\tilde \C_i}|_P=\varphi_P\frac{dz_i}{q_P}=0$
(that is, at the level of forms, and not only at the level of cohomology classes by~(\ref{eq-sec})),
where $z_i$ is a local system of coordinates around $P$ (note that there
is no anti-holomorphic component). Since $C_i$ is irreducible one has $\varphi=pC_i$.
\end{proof}

Analogously to the rational case, consider $P\in \Sing \C$ and three local branches 
$\delta_i$, $\delta_j$, and $\delta_k$ belonging to the global components $\C_i$, $\C_j$,
and $\C_k$ respectively.

\begin{proposition}
\label{prop-zero}
Let us assume that $\psi:=\varphi\frac{\gw}{C_0C_iC_jC_k}$ is trivial in $H^2(S_\C)$ 
for $\varphi=\varphi_iC_i+\varphi_jC_j+\varphi_kC_k$,
where $\varphi$ is a homogeneous polynomial of degree $d_i+d_j+d_k-1$. 
In this case $\varphi=0$.
\end{proposition}

\begin{proof}
Since $\Res^{[2]}\psi=0$, one has that $\psi\in W_1^2$. Furthermore,
${}^1\tilde R^{[1]}_0(\psi)=0$, and hence, by Lemma~\ref{lema-key}, 
$\varphi=fC_iC_jC_k$, which, by the degree condition, implies $f=0$.
\end{proof}

Finally, one can prove that there is a choice of forms satisfying the desired relations.

\begin{theorem}
\label{prop-prod}
There is a choice of representatives $\{\psi_P^{\delta_1,\delta_2}\}_{P,\delta_1,\delta_2}$ in $v_1$ 
and $\{\psi_\infty^{i,k}\}_{i,k}$ in $v_\infty$, such that the following equalities of 2-forms hold:

\begin{equation}
\label{eq-prodijk}
\psi_P^{\delta_1,\delta_2}+\psi_P^{\delta_2,\delta_3}+\psi_P^{\delta_3,\delta_1}=0
\end{equation}
for any $P\in \C_i \cap \C_j \cap \C_k$, $\delta_1$, $\delta_2$, and $\delta_3$ local branches 
of $\C_i$, $\C_j$, and $\C_k$ respectively, and

\begin{equation}
\label{eq-prodij}
\gs_i\wedge \gs_j=
\sum_{{\tiny \begin{matrix} P\in\Sing(\C_{ij})\\ \delta_1\in \Delta_P(\cC_i),\\ \delta_2\in \Delta_P(\cC_j).\end{matrix}}} 
\mu_P(\delta_1,\delta_2) \psi_P^{\delta_1,\delta_2}
+ d_j\sum_{k=2}^{d_i} \psi_\infty^{i,k} - d_i \sum_{k=2}^{d_j} \psi_\infty^{j,k},
\end{equation}
where $\delta_1$ (resp. $\delta_2$) runs over the local branches of $\C_i$ (resp. $\C_j$)
at $P$, and $\mu_P(\delta_1,\delta_2)$ denotes the intersection number.
\end{theorem}

\begin{proof}
For any $P\in \Sing \C$, one can order the set $\Delta_P$ of local branches and denote by 
$\delta_P$ the first local branch under such ordering. Fixing arbitrary representatives 
$\psi_P^{\delta_P,\delta}$ in $v_1$ for any $\delta\in \Delta_P$ one can complete the choices
of representatives as follows:
$\psi_P^{\delta_1,\delta_2}:=\psi_P^{\delta_P,\delta_2} - \psi_P^{\delta_P,\delta_1}$
for any $\delta_1, \delta_2\in \Delta_P$.

As a consequence of that, if $\delta_1$, $\delta_2$, and $\delta_3$ local branches 
of $\C_i$, $\C_j$, and $\C_k$ at $P$ respectively, then one has 
$$\psi_P^{\delta_1,\delta_2}+\psi_P^{\delta_2,\delta_3}+\psi_P^{\delta_3,\delta_1}=
(\psi_P^{\delta_P,\delta_2} - \psi_P^{\delta_P,\delta_1})+
(\psi_P^{\delta_P,\delta_3} - \psi_P^{\delta_P,\delta_2})+
(\psi_P^{\delta_P,\delta_1} - \psi_P^{\delta_P,\delta_3})=0.
$$

In order to prove the second equality we only have freedom on the choice of 
$\psi_P^{\delta_P,\delta}$ and $\psi_\infty^{j,k}$. 

First of all, note that 
$\gs_i\wedge \gs_j=\Jac(C_i,C_j,C_0) \frac{\gw}{C_0C_iC_j}=\left| 
\begin{matrix}
 C_{i,X} & C_{i,Y} & C_{i,Z} \\
 C_{j,X} & C_{j,Y} & C_{j,Z} \\
 C_{0,X} & C_{0,Y} & C_{0,Z} \\
\end{matrix}
\right| \frac{\gw}{C_0C_iC_j}$. 

Denoting by $\psi_{ij}$ the right-hand side of (\ref{eq-prodij}), it is easy to check 
that $\Res^{[2]}\gs_i\wedge \gs_j=\Res^{[2]}\psi_{ij}$ (see~\cite[Theorem 2.47]{ji-tesis}
for details). As a consequence of that, $\gs_i\wedge \gs_j-\psi_{ij}\in K_{ij}\subset W^2_1$. 

By Lemma~\ref{lema-key} and Proposition~\ref{prop-zero}, one only needs to find representatives 
$\psi_P^{\delta_P,\delta}$ verifying ${}^1R^{[1]}_0\left(\gs_i\wedge \gs_j-\psi_{ij}\right)=0$. 

One can proceed as follows: choose arbitrary representatives $\dot{\psi}_P^{\delta_P,\delta}$ 
and $\dot{\psi}_\infty^{i,k}$. We are looking for 2-forms
$\psi_P^{\delta_P,\delta}:=\dot{\psi}_P^{\delta_P,\delta}+\varphi_{P}^{i,\delta}\frac{\omega}{\C_i}$
and 
$\psi_\infty^{i,k}:=\dot{\psi}_\infty^{i,k}+\varphi_{\infty}^{i,k}\frac{\omega}{\C_i}$
for certain $\varphi_{P}^{i,\delta}, \varphi_{P}^{j,\delta}$, and $\varphi_{\infty}^{i,k}$.
satisfying the equations~(\ref{eq-prodij}).

Since the projection of 
$W_0(\log\langle \bar \C\rangle)\rightmap{{}^1R^{[1]}_0} W_0^{[1]}\to 
\oplus \cO_{\bar \C_i}(\bar \C_i)$ is injective (because it is constant when restricted to the 
exceptional divisors), one only needs to make sure that, for a certain choice of representatives, 
such projections are zero.

This gives rise to an affine system of equations on the vector space $\oplus K_{i}$,
where the variables are $\varphi_{P}^{i,\delta}$ and $\varphi_{\infty}^{i,k}$ mentioned above.
All there is left to do is to check that such system is compatible.

Using the relations in~(\ref{eq-prodijk}) one obtains that
\begin{equation}
\label{eq-eij}
E_{ij}\equiv
\gs_i\wedge \gs_j-\dot{\psi}_{ij}=
\sum_{{\tiny \begin{matrix} P\in\Sing(\C_{ij})\\ \delta_1\in \Delta_P(\cC_i),\\ \delta_2\in \Delta_P(\cC_j).\end{matrix}}} 
\mu_P(\delta_1,\delta_2) \left(\varphi_{P}^{j,\delta_2}\frac{\omega}{C_j} - \varphi_{P}^{i,\delta_1} \frac{\omega}{C_i}\right)
+ d_j\sum_{k=2}^{d_i} \varphi_\infty^{i,k} \frac{\omega}{C_i} - d_i \sum_{k=2}^{d_j} \varphi_\infty^{j,k} \frac{\omega}{C_j}.
\end{equation}
where
\begin{equation}
\label{eq-psiij}
\dot{\psi}_{ij}:=
\sum_{{\tiny \begin{matrix} P\in\Sing(\C_{ij})\\ \delta_1\in \Delta_P(\cC_i),\\ \delta_2\in \Delta_P(\cC_j).\end{matrix}}} 
\mu_P(\delta_1,\delta_2) (\dot{\psi}_{P}^{\delta_P,\delta_2} - \dot{\psi}_{P}^{\delta_P,\delta_1})
+ d_j\sum_{k=2}^{d_i} \dot{\psi}_\infty^{i,k} - d_i \sum_{k=2}^{d_j} \dot{\psi}_\infty^{j,k}.
\end{equation}
Therefore, the given system is compatible, if for any linear combination of $\gs_i\wedge \gs_j-\dot{\psi}_{ij}$ 
that maps to zero, then $\gs_i\wedge \gs_j-\dot{\psi}_{ij}=0$.

After a more careful study of~(\ref{eq-eij}), one can prove that a combination of such equations is zero 
on the right-hand side if and only if it can be written as a linear combination of equations of the form
$E_{\bar i,\bar j}+E_{\bar j,\bar k}+E_{\bar k,\bar i}$, where 
$E_{\bar i,\bar j}=\sum_{i_*\in \bar i}\sum_{j_*\in \bar j} E_{i_*,j_*}$, and 
$\C_{\bar i}, \C_{\bar j}, \C_{\bar k}$ form a combinatorial pencil.

All there is left to do is to check that the left-hand side of~(\ref{eq-eij}) is also zero.
By the definition given in~(\ref{eq-psiij}) it is immediate that 
$\dot{\psi}_{\bar i,\bar j}+\dot{\psi}_{\bar j,\bar k}+\dot{\psi}_{\bar k,\bar i}=0$.

In other to check
$\gs_{\bar i}\wedge \gs_{\bar j}+\gs_{\bar j}\wedge \gs_{\bar k}+\gs_{\bar k}\wedge \gs_{\bar i}=0$
one can proceed as follows. By the Combinatorial Noether Theorem~\cite{ji-Marco-noether} we can assume
that $\alpha \C_{\bar i} + \beta \C_{\bar j}=\C_{\bar k}$. In that case,
$$
\begin{matrix}
C_0 C_{\bar i} C_{\bar j} C_{\bar k} \cdot (\gs_{\bar j}\wedge \gs_{\bar k})=
C_0 C_{\bar i} \cdot \Jac(C_0,C_{\bar j},C_{\bar k})\omega=\\
\\
=C_0 C_{\bar i} \cdot \Jac(C_0,C_{\bar j},\alpha \C_{\bar i} + \beta \C_{\bar j})\omega=
-\alpha C_0 C_{\bar i} \cdot \Jac(C_0,C_{\bar i},C_{\bar j})\omega.
\end{matrix}
$$
Analogously
$$
C_0 C_{\bar i} C_{\bar j} C_{\bar k} \cdot (\gs_{\bar k}\wedge \gs_{\bar i})=
-\beta C_0 C_{\bar j} \cdot \Jac(C_0,C_{\bar i},C_{\bar j})\omega.
$$
Hence 
$$
C_0 C_{\bar i} C_{\bar j} C_{\bar k} \cdot 
(\gs_{\bar i}\wedge \gs_{\bar j}+
\gs_{\bar j}\wedge \gs_{\bar k}+
\gs_{\bar k}\wedge \gs_{\bar i})=0
$$
which ends the proof.
\end{proof}

We will denote by $V_1$ (resp. $V_\infty$) the subspace of $W^2\left( \cA^{\log}_{\PP^2}(\C) \right)$ 
generated by the 2-forms in $v_1$ (resp. $v_\infty$) described in Theorem~\ref{prop-prod} and by $\cV_1$ 
(resp. $\cV_\infty$) their projection on $H^2(S_\C)=W^2/d(W^1)$. 
Now we are in the position to give a decomposition of $H^2(S_\C)$.

\begin{corollary}
Under the above conditions 
$$H^2(S_\C)=\cV^2_\cC \oplus \mathcal K_\C\oplus \overline{\mathcal K}_\C.$$
\end{corollary}

\begin{proof}
By Proposition~\ref{prop-ki}, it is enough to check that 
$H^2(S_\C)/(\mathcal K_\C\oplus \overline{\mathcal K}_\C) \cong H^2(S_\C)/\ker \Res^{[2]}$ is isomorphic to
$\cV_\cC^2:=\cV_1\oplus \cV_\infty$. Notice that the residue map $\Res^{[2]}$ is injective on the quotient. Let us 
consider $\psi\in V_1$. By Proposition~\ref{prop-resij} $\left( \Res^{[2]}\psi\right)_{P^i_k}=0$, for any
$i=1,\dots,r$, and any $k=2,\dots,d_i$. On the other hand, by Proposition~\ref{prop-resinfty} it is immediate 
that any 2-form $\psi\in V_\infty$ satisfies that $\sum_{k=1}^{d_i}\left( \Res^{[2]}\psi\right)_{P^i_k}=0$.
Therefore if $\psi\in V_1\cap V_\infty$, 
$\left( \Res^{[2]}\psi\right)_{P^i_k}=0$, $k=1,\dots,d_i$ and thus $\psi=0$.
\end{proof}

As a consequence of the previous results one obtains a description of the cohomology ring 
of the complement of a projective plane curve $H^*(S_\cC)$.

\begin{theorem}
\label{thm-ring}
The cohomology ring $H^*(S_\cC)$ of $S_\cC$ can be decomposed as follows
$$\mathcal V_\cC \oplus \mathcal K_\C\oplus  \overline{\mathcal K}_\C$$
where $H^*(S_\cC)$ is trivial in degree $\geq$ 3, $\mathcal K_\C$ and $ \overline{\mathcal K}_\C$ are 
homogeneous subrings of degree 2 and dimension $g$ each, and $\mathcal V_\cC$ can be described as:
\begin{itemize}
\item Generated in degrees 1 and 2 by
\begin{enumerate}[$(G1)$]
\item (Generators of $\mathcal V_\cC^1$)
$\sigma_1,\dots,\sigma_r,$
\item (Generators of $\mathcal V_\cC^2=\cV_1\oplus \cV_\infty$)
$v_1 \cup v_\infty$ from Theorem~\ref{prop-prod}.
\end{enumerate}
\item The relations given in~(\ref{eq-prodij}), (\ref{eq-prodijk}), and
\begin{equation}
\label{eq-rel1}
\psi_P^{\delta_1,\delta_2}+\psi_P^{\delta_2,\delta_1}=0
\end{equation}
is a complete system of relations
\end{itemize}
$\Box$
\end{theorem}

Note that the ring structure depends on some local and global data which will be described
in what follows. Because of the general condition about the transversal line we will repeat the
Definition~\ref{def-wct} with a slight change in notation.

Let $\tilde \cC=\cC_1\cup ...\cup \cC_r \subset \PP^2$, $\C_0$ a transversal line, $\C:=\C_0 \cup \tilde \C$.
Analogously to Definition~\ref{def-wct} consider the following.

\begin{definition}
\label{def-wct2}
The family $W_\cC:=(\mathbf r,S,\Delta, \phi ,\bar d,\bar g)$ is called the \emph{weak combinatorial type of 
$\cC$ with respect to $\C_0$} or simply \emph{weak combinatorial type of $\cC$} if no ambiguity seems likely to arise.
\end{definition}

\begin{corollary}
\label{cor-ring-wct}
The cohomology ring $H^*(S_\cC)$ of $S_\cC$ only depends on $W_\cC$ the weak combinatorial type of $\cC$. $\Box$
\end{corollary}

\begin{remark}
Corollary~\ref{cor-ring-wct} is also true in the case that the curve does not contain
a transversal line --\,as we have assumed throughout section~\S\ref{sec-cohom-ring}. 
In this case one can add a transversal line and consider 
$\cC=\tilde \cC \cup \cC_0$. The ring $H^*(S_{\tilde \cC})$ fits in the following 
exact sequence 
$$\array{ccccccccc}
0 & \to & H^*(S_{\tilde \cC}) & \to & H^*(S_{\cC}) &
\longrightmap{\pi_{\cC_0}\circ \Res^{[1]}} & \CC_{\cC_0} & \to & 0, 
\endarray$$
where $\Res^{[1]}$ is the residue defined in Proposition~\ref{propinj} 
and $\pi_{\cC_0}$ is the projection of $H^0((\tilde \cC_0\cup \overline{\cC})^{[1]})$ on the 
coordinate corresponding to $\cC_0$.
\end{remark}

\begin{example}
Consider the two conics $\cC_1:=\{y(y-z)+(x+y)^2=0\}$, $\cC_2:=\{y(y-z)+(x-y)^2=0\}$ 
and the line $\cC_3:=\{y=0\}$ (see Figure~\ref{fig-conics}). The weak combinatorial 
type of $\tilde \cC:=\cC_1\cup \cC_2 \cup \cC_3$
is $W_{\tilde \cC}:=(\{1,2,3\},S,\{\Delta_P,\phi_P,(\bullet,\bullet)_P\}_{P\in S},(2,2,1),(0,0,0))$, where
$S:=\{P_1,P_2,P'_2,P_3\}$, $\Delta_{P_1}:=\{\delta^1_1,\delta^1_2\}$,
$\Delta_{P_2}:=\{\delta^2_2,\delta^2_3\}$, 
$\Delta_{P'_2}:=\{\delta'^2_1,\delta'^2_3\}$, $\Delta_{P_2}:=\{\delta^3_1,\delta^3_2\}$,
$\phi_{P_i}(\delta^i_j):=j$, and $(\delta^i_j,\delta^i_k)_{P_i}=i$. 
The ring $H^*(S_{\tilde \cC})$ is generated by the 
1-forms $\omega_i:=2\sigma_3-\sigma_i$, $i=1,2$ and the 2-forms 
$\psi_1:=\psi_{P_3}^{\delta^3_1,\delta^3_2}
+\psi_{P_2}^{\delta^2_2,\delta^2_3}
-\psi_{P'_2}^{\delta'^2_1,\delta'^2_3}$, and
$\psi_2:=\psi_{P_1}^{\delta^1_1,\delta^1_2}
+\psi_{P_2}^{\delta^2_2,\delta^2_3}
-\psi_{P'_2}^{\delta'^2_1,\delta'^2_3}$. The only relation is given by
$\omega_1\wedge \omega_2=3\psi_1+\psi_2$. Hence
$$
H^*(S_{\tilde \cC})=\langle \omega_1, \omega_2, \psi_1, \psi_2 \mid
\omega_1\wedge \omega_2=3\psi_1+\psi_2
\rangle.
$$
\begin{figure}[ht]
\includegraphics[scale=.5]{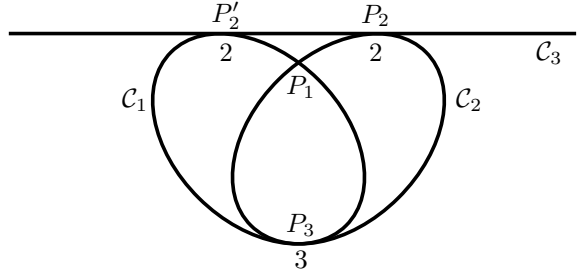}
\begin{picture}(0,0)
\put(-170,55){\makebox(0,0){$\cC_1$}}
\put(-45,55){\makebox(0,0){$\cC_2$}}
\put(-15,73){\makebox(0,0){$\cC_3$}}
\put(-108,8){\makebox(0,0){$P_3$}}
\put(-108,-5){\makebox(0,0){$3$}}
\put(-108,60){\makebox(0,0){$P_1$}}
\put(-80,73){\makebox(0,0){$2$}}
\put(-80,87){\makebox(0,0){$P_2$}}
\put(-136,73){\makebox(0,0){$2$}}
\put(-136,87){\makebox(0,0){$P'_2$}}
\end{picture}
\caption{Projective realization of $\cC$ and multiplicities of intersection}
\label{fig-conics}
\end{figure}
\end{example}

\begin{remark}
\label{rem-affine-proj}
As in Definition~\ref{def-wct2}, the curve $\tilde \C\subset \PP^2$ will not be assumed to
have a transversal line and usually, we will denote by $\C$ the union of $\tilde \C$ and
a transversal line. In the future, we will always consider this situation unless otherwise stated.
\end{remark}

\section{Formality of complements to projective plane curves}

All the basic definitions of minimal algebras, minimal models, homotopy, etc required in the definition of 
formality and in the theory of homotopy theory of algebras can be found for instance in any of the foundational 
papers~\cite{Deligne-Griffiths-Morgan-Sullivan-real-homotopy,Morgan-algebraic-topology,Sullivan-infinitesimal}.

\begin{definition}
Two graded differential algebras $(A,d_A)$ and $(B,d_B)$ are called \emph{quasi-isomorphic} if there exists a 
morphism of graded algebras $f:A\to B$ such that the induced morphism $f^*:H^*(A,d_A)\to H^*(B,d_B)$ is an isomorphism.
Note that ``being quasi-isomorphic" is not an equivalence relation. We will refer to the \emph{quasi-isomorphism class}
of a graded differential algebra as the minimal equivalence class generated by the quasi-isomorphism relation.

A minimal differential graded algebra is called \emph{formal} if it is quasi-isomorphic to its cohomology algebra 
using a zero differential. A differential graded algebra is called \emph{formal} if its minimal model is formal. 
Finally, a complex space $X$ is called \emph{formal} if the algebra of differential forms $(\cE(X),d)$ is formal.
\end{definition}

The concept of formal algebra is well defined since any differential graded algebra has a unique 
(up to homotopy) minimal model (c.f.~\cite[Section \S 5]{Sullivan-infinitesimal}). Also note that a minimal model 
for $(A,d_A)$ consists of a minimal algebra $(\cM(A),d_{\cM(A)})$ plus a quasi-isomorphism $\cM(A)\to A$. Therefore, 
if one finds a quasi-isomorphism between $(\cE(X),d)$ and $(H(X),0)$ then $X$ is formal. Moreover, if $X$ is a 
smooth complex variety and $\overline X$ is a completion of $X$ by a simple normal crossing divisor, then the minimal model 
of $\cE(X)$ and the minimal model of $\cA_{\overline X}(\log\langle D \rangle)$ are isomorphic 
(c.f.~\cite[Section \S 6]{Morgan-algebraic-topology}).

\begin{theorem}
There is a well-defined quasi-isomorphism 
$H^*(S_{\cC}) \to \cA^*_{\overline S_\cC}(\log\langle \overline \cC \rangle)$.
\end{theorem}

\begin{proof}
According to Theorem~\ref{thm-ring} $H^*(S_{\cC})$ admits a decomposition 
$$
\CC\oplus \cV_\cC^1 \oplus \cV_\cC^2 \oplus \cK_\cC \oplus \overline \cK_\cC,
$$
where $\cV_\cC^1$ is generated by 1-forms $V_\cC^1\subset W^1$, and
$\cV_\cC^2$ (resp. $\cK_\cC$, $\overline \cK_\cC$) are generated by 2-forms $V_\cC^2\subset W_2^2$
(resp. $K_\cC\subset W^2_1$, $\overline K_\cC\subset W^2_1$).
Each cohomology class $\varphi\in H^*(S_\cC)$ can be thus described as the cohomology class of 
a combination of forms as follows: $\varphi = [z+ \psi^1 + \psi^2_1+ \psi^2_2]$. The map is defined
by $\varphi \mapsto z+ \psi^1 + \psi^2_1+ \psi^2_2$. By Theorem~\ref{prop-prod} this map is well
defined and it is obviously a quasi-isomorphism.
\end{proof}

As a consequence of the discussion at the beginning of this section one has the following.

\begin{theorem}
\label{thm-formal}
The complement of a plane projective curve $S_\cC$ is a formal space.
\end{theorem}

\begin{remark}
Theorem~\ref{thm-formal} is the global version of the formality of algebraic links proved by Durfee-Hain
in~\cite{Durfee-mixed-hodge-homotopy}. The result is a consequence of a more general fact proved
paper by A.~Macinic in~\cite{Macinic-cohomology}: a $2$-complex $X$ which is $1$-formal is 
also a formal space.

The $1$-minimal model $\mathcal{M}_1(A)$ of a differential graded algebra $(A,d_A)$ is the 
subalgebra generated by the degree $1$ part in $\mathcal{M}(A)$. Then a space $X$ is 
$1$-formal if $\mathcal{M}_1(\mathcal{E}(X),d)$ is quasi-isomorphic to $\mathcal{M}_1(H^*(X),0)$.
This condition can be restated in terms of the fundamental group as follows.
A finitely presented group $G$ is $1$-formal if its Malcev completion is filtered isomorphic 
to its holonomy Lie algebra, completed  with respect to bracket length.
Fundamental groups of complements to algebraic plane curves are known to be $1$-formal, 
(see~\cite{Kohno-holonomyLie} and \cite{Morgan-algebraic-topology}).
\end{remark}

\section{Examples}
\label{sec-ex}

\subsection{Weak combinatorics does not determine classical combinatorics}
\label{ex-wc-no-kc}
\begin{figure}[ht]
\includegraphics[scale=.5]{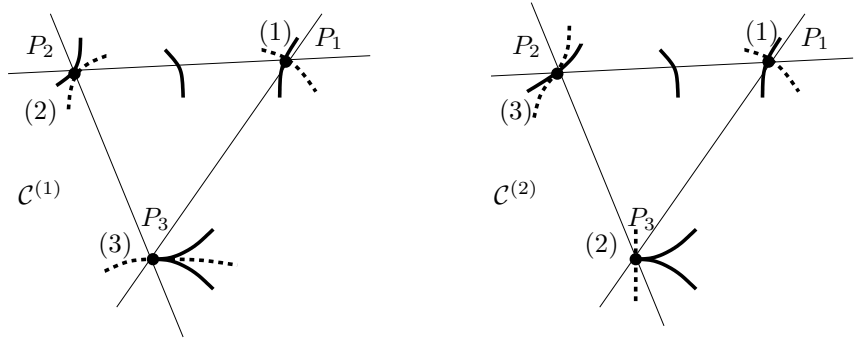}
\begin{picture}(0,0)
\put(-308,110){\makebox(0,0){$P_2$}}
\put(-308,85){\makebox(0,0){$(2)$}}
\put(-200,113){\makebox(0,0){$P_1$}}
\put(-220,115){\makebox(0,0){$(1)$}}
\put(-265,45){\makebox(0,0){$P_3$}}
\put(-280,35){\makebox(0,0){$(3)$}}
\put(-308,55){\makebox(0,0){$\C^{(1)}$}}

\put(-126,110){\makebox(0,0){$P_2$}}
\put(-130,85){\makebox(0,0){$(3)$}}
\put(-18,113){\makebox(0,0){$P_1$}}
\put(-38,115){\makebox(0,0){$(1)$}}
\put(-83,45){\makebox(0,0){$P_3$}}
\put(-98,35){\makebox(0,0){$(2)$}}
\put(-130,55){\makebox(0,0){$\C^{(2)}$}}
\end{picture}
\caption{Singularities of $\C^{(1)}$ and $\C^{(2)}$ respectively.}
\label{fig-wc-kc}
\end{figure}

Let $\ell_0:=\{x=0\}$, $\ell_1:=\{y=0\}$, and $\ell_2:=\{z=0\}$ be three lines in general
position and consider: 
\begin{enumerate}
\item
$\tilde \C_1:=\{(x-y)^2-(x+y)z=0\}$ a conic tangent to $\ell_2$ at 
$(1,1,0)$ and passing through $\ell_0\cap \ell_1$, 
\item
$\tilde \C^{(1)}_2:=\{x-y+z=0\}$ the line
passing through $\ell_0\cap \tilde \C_1$ and $\ell_2\cap \tilde \C_1$, and 
\item
$\tilde \C^{(2)}_2:=\{3x-y+z=0\}$ the line tangent to $\tilde \C_1$ at 
$\ell_0\cap \tilde \C_1$. 
\end{enumerate}
The Cremona transformation based on $\ell_0$, $\ell_1$, and $\ell_2$
transforms $\tilde \C^{(k)}:=\tilde \C_1\cup \tilde \C^{(k)}_2$ into $\C^{(k)}$, a union
of a cuspidal cubic $\C_1$ and a conic $\C^{(k)}_2$. 
Note that $\C^{(k)}$ has three singular points $\{P_1,P_2,P_3\}$, where 
$\Delta_{\C^{(k)},P_i}:=\{\delta^1_i,\delta^{2,k}_i\}$, 
$\phi_{\C^{(k)},P_i}(\delta^{1}_i)=\C_1$, 
$\phi_{\C^{(k)},P_i}(\delta^{2,k}_i)=\C^{(k)}_2$, 
$(\delta^1_i,\delta^{2,1}_i)_{\C^{(1)},P_i}=i$,
and $(\delta^1_i,\delta^{2,2}_i)_{\C^{(2)},P_i}=\sigma_{(2,3)}(i)$ 
(where $\sigma_{(2,3)}(i)$ represents the permutation $(2,3)$ applied to $i$).
Figure~\ref{fig-wc-kc} represents the singular points of $\C^{(k)}$, the local branches 
at those points (solid line for $\C_1$ and broken line for $\C_2^{(k)}$), and the 
multiplicity of intersection in brackets.
Note that the bijection $\varphi$ of singular points that permutes $P_2$ and $P_3$ induces an
equivalence between $W_{\C^{(1)}}$ and $W_{\C^{(2)}}$, since
$\varphi_{\mathbf r_{\C^{(1)}}}$ and $\varphi_{P_i}$ are forced by their compatibility with 
the degrees and with $\varphi$. The combinatorial types $K_{\C^{(1)}}$ and $K_{\C^{(2)}}$ 
cannot be equivalent since the topological types of their singularities do not coincide.

\subsection{An explicit computation of the cohomology ring in the non-rational case}
\begin{figure}[ht]
\includegraphics[scale=.5]{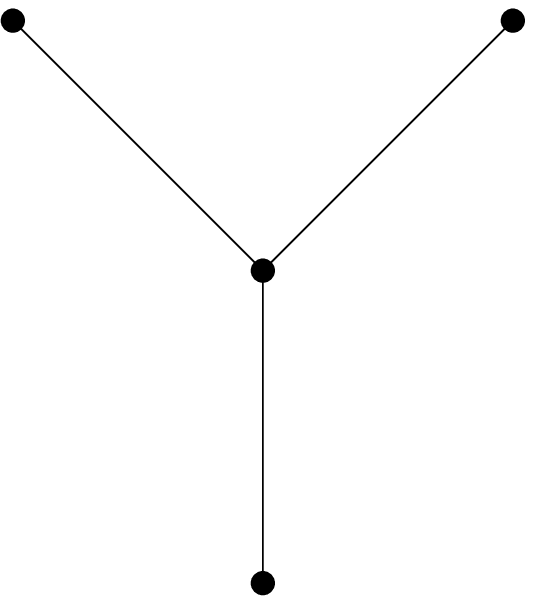}
\begin{picture}(0,0)
\put(-30,-5){\makebox(0,0){$1,1$}}
\put(-30,40){\makebox(0,0){$1,0$}}
\put(-95,80){\makebox(0,0){$0,0$}}
\put(12,80){\makebox(0,0){$0,0$}}
\end{picture}
\caption{Description of $\cT^{\nul}_{\AAA_3}$, and $\cT_{\AAA_3}$ respectively.}
\label{fig-tacnode}
\end{figure}

We will present a simple example of a non-rational arrangement of curves in order to
show how to compute the forms described in \S\ref{sec-cohom-ring}.
Let $\cC:=\cC_0\cup \cC_1 \cup  \cC_2 \cup \cC_3$, where 
$\cC_0:=\{x+y+z=0\}$, $\cC_1:=\{y-z=0\}$, $\cC_2:=\{xy+xz+yz=0\}$, 
and $\cC_3:=\{x^2(y+z)+y^2(x+z)+z^2(x+y)=0\}$.
In this case, for simplicity it is more
convenient to consider the line at infinity $\cC_0$ with an equation different from $\{z=0\}$. 
Consider $\xi$ a primitive third root of unity (a root of $t^2+t+1=0$) and denote 
$\cC_0\cap \cC_1=\{P_{01}=[-2:1:1]\}$,
$\cC_0\cap \cC_2=\{R_1=[-\bar \xi-1:\bar \xi:1],R_2=[-1-\xi:\xi:1]\}$,
$\cC_0\cap \cC_3=\{Q_1=[0:1:-1],Q_2=[-1:0:1],Q_3=[-1:1:0]\}$,
$\cC_1\cap \cC_2=\{P_1,P_{12}=[1:-2:-2]\}$,
$\cC_1\cap \cC_3=\{P_1,P_{13}=[\xi:1:1],\bar P_{13}=[\bar \xi:1:1]\}$,
$\cC_2\cap \cC_3=\{P_1=[1:0:0],P_2=[0:1:0],P_3=[0:0:1]\}$,

Since all the local branches of the irreducible components at any singular point are irreducible, we will denote
by $\psi_P^{i,j}$ the 2-form associated with the singular point $P$ and the local branches at $P$ of $\cC_i$ and $\cC_j$. 
For example, in order to compute $\psi_{P_1}^{2,3}=\varphi_{P_1}^{2,3}\frac{\omega}{C_0C_2C_3}$, one needs a section
of $H^0(\PP^2,\cI^{2,3}_{P_1}(3))$. Note that
$$
\left(\cI^{2,3}_{P_1}\right)_{P}:=
\begin{cases}
\{\varphi\in \cO_{P} \mid \cT_P(\C_{2,3})|_\varphi \geq \cT_{\AAA_3}\}=\mathfrak m_P & \text{ if } P=P_1,\\
\{\varphi\in \cO_{P} \mid \cT_P(\C_{2,3})|_\varphi \geq \cT_{\AAA_3}^{\nul}\} & \text{ if } P=P_2,P_3,\\
\cO_{P} & \text{ otherwise, }
\end{cases}
$$
where Figure~\ref{fig-tacnode} describes the local conditions at the tacnodes $P_i$.

Therefore $\varphi_{P_1}^{2,3}$ is the equation of a cubic 
$\alpha(xz+x^2+(1-\xi)xy+yz)z+\beta C_0C_2$. 
In order to obtain a normal form one has to require the different residues of $\psi_{P_1}^{2,3}$ 
at $P_1$ and at an exceptional divisor $E$ joining $\delta_2$ and $\delta_3$ to be equal to $\pm 1$. 
It is a simple computation that
$$\Res^{[2]}_{P_1}\psi_{P_1}^{2,3}=\frac{\alpha}{3}$$ 
and that
$$\left({}^1R_0^{[1]}\psi_{P_1}^{2,3}\right)_{E}=\frac{1}{3}(\beta-\xi).$$ 
Since 
$\left({}^1R_0^{[1]} \frac{d\cC_2}{\cC_2}\wedge \frac{d\cC_3}{\cC_3}\right)_{E}=-\frac{2}{3}$
and $(\delta_2,\delta_3)_{P_1}=2$, one concludes that
$$\varphi_{P_1}^{2,3}=3(xz+x^2+(1-\xi)xy+yz)z+(\xi-1) C_0C_2.$$

Analogously one can proceed with $\psi_{P_1}^{1,2}=\varphi_P^{1,2}\frac{\omega}{C_0C_1C_2}$ and
$\psi_{P_1}^{3,1}=\varphi_{P_1}^{3,1}\frac{\omega}{C_0C_3C_1}$ obtaining $\varphi_{P_1}^{1,2}:=2x-\xi y+(1+\xi)z$ 
and $\varphi_{P_1}^{1,3}:=2(x^2+xz+2yz-y^2)+C_0C_1$.

Note that $\varphi_{P_1}^{1,2}C_3+\varphi_{P_1}^{2,3}C_1-\varphi_{P_1}^{1,3}C_2=0$ and hence
$\psi_{P_1}^{1,2}+\psi_{P_1}^{2,3}+\psi_{P_1}^{3,1}=0$ (Theorem~\ref{prop-prod}(\ref{eq-prodijk})).

The following list describes the polynomials $\varphi_P^{i,j}$ for the generating 2-forms 
$\psi_P^{i,j}:=\varphi_P^{i,j}\frac{\omega}{C_0C_iC_j}$:
$$
\begin{array}{l}
\varphi_{P_1}^{2,3}:=(\xi+2)(zx^2+\xi yx^2+xz^2+\xi y^2x+z^2y+\xi zy^2),\\
\varphi_{P_2}^{2,3}:=
(\xi+2)(y^2x+xz^2+yx^2+zx^2+z^2y+ (1-\xi) zxy+\xi zy^2),\\
\varphi_{P_3}^{2,3}:=
(\xi-1)(y^2x+xz^2+yx^2+zx^2+y^2z+(\xi + 2)zxy- (1+\xi) z^2y),\\
\varphi_{P_1}^{1,2}:=2x-\xi y+(1+\xi)z,\\
\varphi_{P_{12}}^{1,2}:=(\xi-1)((\xi+1)y+z),\\
\varphi_{P_{1}}^{1,3}:=2x^2+xz+xy-z^2+4yz-y^2,\\
\varphi_{P_{13}}^{1,3}:=
-(\xi + 2)(xz+xy+z^2+2\xi zy+y^2),\\
\varphi_{\bar P_{13}}^{1,3}:=
(\xi-1)(xz+xy+z^2-2(\xi + 1)yz+y^2).\\
\end{array}
$$
Finally we also describe the polynomials $\varphi_\infty^{i,k}$ for the generating 2-forms 
$\psi_\infty^{i,k}:=\varphi_\infty^{i,k}\frac{\omega}{C_0C_i}$:
$$
\begin{array}{l}
\varphi_\infty^{3,Q_2}:=
-3(x+y),\\
\varphi_\infty^{3,Q_3}:=
3(x+z),\\
\varphi_\infty^{2,R_2}:=-(2\xi+1).\\
\end{array}
$$
One can then easily verify that
$$
\Jac(C_2,C_3,C_0)=2\varphi_{P_1}^{2,3}+2\varphi_{P_2}^{2,3}+2\varphi_{P_3}^{2,3}+
3\varphi_\infty^{2,R_2}C_3-2\varphi_\infty^{3,Q_2}C_2-2\varphi_\infty^{3,Q_3}C_2,
$$
and hence
$$
\sigma_{2}\wedge \sigma_{3}=2\psi_{P_1}^{2,3}+2\psi_{P_2}^{2,3}+2\psi_{P_3}^{2,3}+
3\psi_\infty^{2,R_2}-2\psi_\infty^{3,Q_2}-2\psi_\infty^{3,Q_3},
$$
which corresponds to Theorem~\ref{prop-prod}(\ref{eq-prodij}).

The same can be checked for $\Jac(C_1,C_2,C_0)$ and $\Jac(C_1,C_3,C_0)$.

\bibliographystyle{amsplain}
\def\cprime{$'$}
\providecommand{\bysame}{\leavevmode\hbox to3em{\hrulefill}\thinspace}
\providecommand{\MR}{\relax\ifhmode\unskip\space\fi MR }
\providecommand{\MRhref}[2]{%
  \href{http://www.ams.org/mathscinet-getitem?mr=#1}{#2}
}
\providecommand{\href}[2]{#2}

\end{document}